\ifdefined\singlecolumn
\documentclass[11pt,draftclsnofoot,onecolumn]{IEEEtran}
\else
\documentclass[10pt,twocolumn]{IEEEtran}
\fi

\usepackage{graphicx,cite}
\usepackage[tbtags]{amsmath}
\usepackage{amssymb,color,amsthm}
\usepackage[none]{hyphenat}
\usepackage{subfigure}
\usepackage{bbm}
\usepackage{etoolbox}
\usepackage[dvipsnames]{xcolor}
\usepackage{url}

\newtoggle{twocolumn}
\ifCLASSOPTIONtwocolumn\toggletrue{twocolumn}\else\togglefalse{twocolumn}\fi
\def\onlytwo#1{\iftoggle{twocolumn}{#1}{}}
\def\onlyone#1{\iftoggle{twocolumn}{}{#1}}

\def\bx{\mathbf x}

\def\bw{\mathbf w}
\def\by{\mathbf y}
\def\byt{\widehat{\mathbf y}}

\def\bA{\mathbf A}

\def\bI{\mathbf I}

\def\bM{\mathbf M}

\def\bP{\mathbf P}

\def\bU{\mathbf U}

\def\bW{\mathbf W}
\def\bX{\mathbf X}
\def\bY{\mathbf Y}
\def\cA{\mathcal A}
\def\hcA{\widehat{\mathcal A}}
\def\cB{\mathcal B}
\def\cC{\mathcal C}
\def\bYt{\widehat{\mathbf Y}}
\def\bLambda{\mathbf \Lambda}
\def\bOmega{\mathbf \Omega}
\def\bZ{\mathbf Z}
\newcommand{\T}{\scriptscriptstyle T}

\def\leref#1{Lemma~\ref{#1}}
\def\thref#1{Theorem~\ref{#1}}

\def\figref#1{Figure~\ref{#1}}
\def\figtab#1{Table~\ref{#1}}
\def\bydef{\triangleq}
\newtheorem{lemma}{Lemma}
\newtheorem{theorem}{Theorem}
\newtheorem{corollary}{Corollary}

\begin{document}

\title{A Nonconvex Splitting Method for Symmetric Nonnegative Matrix Factorization: \\Convergence Analysis and Optimality}
\author{Songtao~Lu, \emph{Student Member, IEEE},
    Mingyi~Hong, \emph{Member, IEEE}
    and Zhengdao~Wang, \emph{Fellow, IEEE}

\thanks{Manuscript received May 15, 2016; revised October 6, 2016,
January 6, 2017, and February 16, 2017; accepted February 20, 2017. The associate editor coordinating the review of this
manuscript and approving it for publication was Marco Moretti. Part of the
paper was presented at the 42nd IEEE International Conference on Acoustics, Speech, and Signal Processing (ICASSP), New Orleans, March 5--9, 2017. This work was
supported in part by NSF under Grants No.~1523374 and No.~1526078, and by
AFOSR under Grant No.~15RT0767.}

\thanks{
    Songtao Lu and Zhengdao Wang are with the Department of Electrical and Computer Engineering,
    Iowa State University, Ames, IA 50011, USA (emails: \{songtao, zhengdao\}@iastate.edu).}

\thanks{
    Mingyi Hong is with the Department of Industrial and Manufacturing Systems Engineering,
    Iowa State University, Ames, IA 50011, USA (email: mingyi@iastate.edu). }
}

\maketitle

\begin{abstract}
Symmetric nonnegative matrix factorization (SymNMF) has important
applications in data analytics problems such as document clustering, community
detection and image segmentation. In this paper, we propose a novel nonconvex
variable splitting method for solving SymNMF. The proposed algorithm is
guaranteed to converge to the set of Karush-Kuhn-Tucker (KKT) points of the
nonconvex SymNMF problem. Furthermore, it achieves a global sublinear
convergence rate. We also show that the algorithm can be efficiently
implemented in parallel. Further, sufficient conditions are provided which
guarantee the global and local optimality of the obtained solutions. Extensive
numerical results performed on both synthetic and real data sets suggest that
the proposed algorithm converges quickly to a local minimum solution.

\begin{keywords}
Symmetric nonnegative matrix factorization, Karush-Kuhn-Tucker points,
variable splitting, global and local optimality, clustering
\end{keywords}
\end{abstract}

\section{Introduction} \label{sec:intro}

Nonnegative matrix factorization (NMF) refers to factoring a given matrix into
the product of two matrices whose entries are all nonnegative. It has long
been recognized as an important matrix decomposition problem
\cite{capo81,pata94}. The requirement that the factors are component-wise
nonnegative makes NMF distinct from traditional methods such as the principal
component analysis (PCA) and the linear discriminant analysis (LDA), leading
to many interesting applications in imaging, signal processing and machine
learning \cite{giva14,wazh13,hoye04,lese01,yafu15}; see
\cite{Gillis15} for a recent survey. When further requiring that the two
factors are identical after transposition, NMF becomes the so-called symmetric
nonnegative matrix factorization (SymNMF). In the case where the given matrix
cannot be factorized exactly, an approximate solution with a suitably defined
approximation error is desired. Mathematically, SymNMF approximates a given
(usually symmetric) nonnegative matrix $\bZ\in\mathbb{R}^{N\times N}$ by a low
rank matrix $\bX\bX^{\T}$, where the factor matrix
$\mathbf{X}\in\mathbb{R}^{N\times K}$ is component-wise nonnegative, typically
with $K\ll N$. Let $\|\cdot\|_F$ denote the Frobenius norm. The problem can be
formulated as a nonconvex optimization problem \cite{zhshhe12,kesi14,kuyu15}:
\begin{equation}\label{eq:symnmf}
\min_{\bX\ge 0} \quad f(\bX)=\frac{1}{2}\|\bX\bX^{\T}-\bZ\|^2_F.
\end{equation}

Recently, SymNMF has found many applications in document clustering, community
detection, image segmentation and pattern clustering in bioinformatics
\cite{kuyu15,wali11,zhshhe12}. An important class of clustering methods is
known as spectral clustering, e.g., \cite{luxb07,luwa15}, which is based on
the eigenvalue decomposition of some transformed graph Laplacian matrix. In
\cite{ding05}, it has been shown that spectral clustering and SymNMF are two
different ways of relaxing the kernel $K$-means clustering, where the former
relaxes the nonnegativity constraint while the latter relaxes certain
orthogonality constraint. SymNMF also has the advantage of often yielding more
meaningful and interpretable results~\cite{kuyu15}.

\subsection{Related Work}

Due to the importance of the NMF problem, many algorithms have been proposed
in the literature for finding its high-quality solutions. Well-known
algorithms include the multiplicative update \cite{lese01}, alternating
projected gradient methods \cite{lich07}, alternating nonnegative least
squares (ANLS) with the active set method \cite{kipa11} and a few recent
methods such as the bilinear generalized approximate message passing
\cite{pasc14i,pasc14ii}, as well as methods based on the block coordinate
descent \cite{kihe13}. These methods often possess strong convergence
guarantees (to Karush-Kuhn-Tucker (KKT) points of the NMF problem) and most of
them lead to satisfactory performance in practice; see \cite{Gillis15} and the
references therein for detailed comparison and comments for different
algorithms. Unfortunately, most of the aforementioned methods for NMF lack
effective mechanisms to enforce the symmetry between the resulting factors,
therefore they are not directly applicable to SymNMF. Recently, there have
been works focusing on customized algorithms for SymNMF, which we
review below.

To this end, first rewrite SymNMF equivalently as
\begin{equation}\label{eq:symnmf:2}
\min_{\bY\ge 0,\; \bX=\bY} \frac{1}{2}\|\bX\bY^{\T}-\bZ\|^2_F.
\end{equation}
A simple strategy is to ignore the equality constraint $\bX=\bY$, and then
alternatingly perform the following two steps: 1) solving $\bY$ with $\bX$
being fixed (a nonnegative least squares problem); 2) solving $\bX$ with $\bY$
being fixed (a least squares problem). Such ANLS algorithm has been proposed
in \cite{kuyu15} for dealing with SymNMF. Unfortunately, despite the fact that
an optimal solution can be obtained in each subproblem, there is no guarantee
that the $\bY$-iterate will converge to the $\bX$-iterate. The algorithm in
\cite{kuyu15} adds a regularized term for the difference between the two
factors to the objective function and explicitly enforces that the two
matrices are equal at the output. Such an extra step enforces symmetry, but
unfortunately also leads to the loss of global convergence guarantees. A
related ANLS-based method has been introduced in
\cite{kesi14}; however the algorithm is based on the assumption that there
exists an exact symmetric factorization (i.e., $\exists~\bX\ge 0$ such that
$\bX\bX^{\T}=\bZ$). Without such assumption, the algorithm may not converge to
the set of KKT points\footnote{\label{f.con}Let $d(a,s)$ denote the distance
between two points $a$ and $s$. We say that a sequence $a_i$ converges to a
set $\mathcal{S}$ if the distance between $a_i$ and $\mathcal{S}$, defined as
$\inf_{s\in \mathcal{S}} d(a_i, s)$, converges to zero, as $i\to\infty$.} of problem
\eqref{eq:symnmf}. A multiplicative update for SymNMF has been proposed in
\cite{zhshhe12}, but the algorithm lacks convergence guarantees (to KKT points
of problem \eqref{eq:symnmf})
\cite{lin07}, and has a much slower convergence speed than the one proposed in
\cite{kesi14}. In
\cite{kuyu15,kudi12}, algorithms based on the projected gradient descent (PGD)
and the projected Newton (PNewton) have been proposed, both of which directly
solve the original formulation \eqref{eq:symnmf}. Again there has been no
global convergence analysis since the objective function is a nonconvex
fourth-order polynomial. More recently, the work \cite{vagi16} applies the
nonconvex coordinate descent (CD) algorithm for SymNMF. Due to the fact that
the minimizer of the fourth order polynomial is not unique in each coordinate
updating, the CD-based method may not converge to stationary points.

Another popular method for NMF is based on the alternating direction method of
multipliers (ADMM), which is a flexible tool for large scale convex
optimization \cite{BoydADMMsurvey2011}. For example, using ADMM for both NMF
and matrix completion, high quality results have been obtained in
\cite{xuyi12} for gray-scale and hyperspectral image recovery. Furthermore,
ADMM has been applied to generalized versions of NMF where the objective
function is the general beta-divergence \cite{sufe14}. A hybrid alternating
optimization and ADMM method was proposed for NMF, as well as tensor
factorization, under a variety of constraints and loss measures in
\cite{husi15}. However, despite the promising numerical results, none of the
works discussed above has rigorous theoretical justification for SymNMF.
Recently, the work \cite{hach16} has applied the ADMM for NMF and provided one
of the first analysis for using ADMM to solve nonconvex matrix-factorization
type problems. However, it is important to note that the algorithm in
\cite{hach16} does not apply to the SymNMF case, because our problem is more
restrictive in that symmetric factors are desired, while in NMF symmetry is
not enforced. Technically, imposing symmetry poses much difficulty in the
analysis (we will comment on this point shortly). In fact, the convergence of
ADMM for SymNMF is still open in the literature.

An important research question for NMF and SymNMF is whether it is possible to
design algorithms that lead to {\it globally} optimal solutions. At the first
sight such problem appears very challenging since finding the exact NMF is
NP-hard \cite{vava09} and checking whether a positive semidefinite matrix can
be decomposed exactly by SymNMF is also NP-hard \cite{pelu14}. However, some
promising recent findings suggest that when the structure of the underlying
factors are appropriately utilized, it is possible to obtain rather strong
results. For example, in \cite{sare15}, the authors have shown that for the
low rank factorized stochastic optimization problem where the two low rank
matrices are symmetric, a modified stochastic gradient descent algorithm is
capable of converging to a global optimum with constant probability from a
random starting point. Related works also include
\cite{gillis12, sunluo14,zhao15}. However, when the factors are required to be
nonnegative and symmetric, it is no longer clear whether the existing analysis
can still be used to show convergence to global/local optimal points. For the
nonnegative principal component problem (i.e., finding the leading nonnegative
eigenvector) under the spiked model, reference
\cite{mori15} shows that certain approximate message passing algorithm is able
to find the global optimal solution asymptotically. Unfortunately, this
analysis does not generalize to an arbitrary symmetric observation matrix for
the case $K>1$. To our best knowledge, a characterization of global and local
optimal solutions for SymNMF is still lacking.

\subsection{Contributions}

In this paper, we first propose a novel algorithm for SymNMF, which utilizes
nonconvex splitting and is capable of converging to the set of KKT points with
a provable global convergence rate. The main idea is to relax the symmetry
requirement at the beginning and gradually enforce it as the algorithm
proceeds. Second, we provide a number of easy-to-check sufficient conditions
guaranteeing the local or global optimality of the obtained solutions.
Numerical results on both synthetic and real data show that the proposed
algorithm achieves fast and stable convergence (often to local minimum
solutions) with low computational complexity.

More specifically, the main contributions of this paper are:

1) We design a novel nonconvex splitting SymNMF (NS-SymNMF) algorithm, which
converges to the set of KKT points of SymNMF with a global sublinear rate. To
our best knowledge, it is the first SymNMF solver that possesses global
convergence rate guarantees.

2) We provide a set of easily checkable sufficient conditions (which only
involve finding the smallest eigenvalue of certain matrix) that characterize
the global and local optimality of the solutions. By utilizing such
conditions, we demonstrate numerically that with high probability, our
proposed algorithm converges not only to the set of KKT points but to a local
optimal solution as well.

{\bf Notation:} Bold upper case letters without subscripts (e.g., $\bX,\bY$)
denote matrices and bold lower case letters without subscripts (e.g.,
$\bx,\by$) represent vectors. The notation $\bZ_{i,j}$ denotes the $(i,j)$-th
entry of matrix $\bZ$. Vector $\bX_i$ denotes the $i$th row of matrix $\bX$
and $\bX'_m$ denotes the $m$th column of the matrix. 

\section{The Proposed Algorithm}

The proposed algorithm leverages the reformulation \eqref{eq:symnmf:2}. Our
main idea is to gradually tighten the difficult equality constraint $\bX=\bY$
as the algorithm proceeds so that when convergence is approached, such
equality is eventually satisfied. To this end, let us construct the augmented
Lagrangian for \eqref{eq:symnmf:2}, given by
\begin{equation}\label{eq:aug}
\mathcal{L}(\bX,\bY; \bLambda)=\frac{1}{2}\|\bX\bY^{\T}-\bZ\|^2_F+\langle\bY-\bX,\bLambda\rangle+\frac{\rho}{2}\|\bY-\bX\|^2_F
\end{equation}
where $\bLambda\in\mathbb{R}^{N\times K}$ is a matrix of dual variables,
$\langle\cdot\rangle$ denotes the inner product operator, and $\rho>0$ is a
penalty parameter whose value will be determined later.

It may be tempting to directly apply the well-known ADMM method to the
augmented Lagrangian \eqref{eq:aug}, which alternatingly minimizes the primal
variables $\bX$ and $\bY$, followed by a dual ascent step $\bLambda\leftarrow
\bLambda + \rho (\bY-\bX)$. Unfortunately, the classical result for ADMM
presented in
\cite{BoydADMMsurvey2011,Bertsekas87,EcksteinBertsekas1992} only works for
convex problems, hence they do not apply to our nonconvex problem
\eqref{eq:symnmf:2} (note this is a linearly constrained {\it nonconvex}
problem where the nonconvexity arises in the objective function). Recent
results such as
\cite{hong14nonconvex_admm,li14nonconvex,Ames13LDA,wang15nonconvexadmm} that
analyze ADMM for nonconvex problems do not apply either, because in these
works the basic requirements are: 1) the objective function is separable over
the block variables; 2) the smooth part of the augmented Lagrangian has
Lipschitz continuous gradient with respect to all variable blocks.
Unfortunately neither of these conditions are satisfied in our problem.

Next we begin presenting the proposed algorithm. We start by considering the
following reformulation of problem \eqref{eq:symnmf}
\begin{align}\label{eq:symnmf:3}
\min_{\bX,\bY}&\quad \frac{1}{2}\|\bX\bY^{\T}-\bZ\|^2_F\\
\mbox{s.t.}&\quad \bY\ge 0, \; \bX=\bY, \; \|\bY_i\|^2_2\le \tau, \; \forall~i, \nonumber
\end{align}
where $\tau>0$ is some given constant.

Let $\bOmega^*$ denote the dual matrix for the constraint
$\bX\ge 0$ in the Lagrangian of problem \eqref{eq:symnmf}. The KKT conditions
of problem \eqref{eq:symnmf} are given by \cite[eq. (5.49)]{bova04}
\begin{subequations}
\begin{align}
&2\left(\bX^*(\bX^*)^{\T}-\frac{\bZ^{\T}+\bZ}{2}\right)\bX^*-\bOmega^*=0, \label{eq.kktgrad}
\\
&\bOmega^*\ge 0,\label{eq.kktopto}
\\
&\bX^*\ge 0, \label{eq.kktoptx}
\\
&\bX^*\circ\bOmega^*=\mathbf{0} \label{eq.kktslack}
\end{align}
\end{subequations}
where $\circ$ denotes the Hadamard product. For a point $\bX^*$, if we can
find some $\bOmega^*$ such that $(\bX^*, \bOmega^*)$ satisfies conditions
\eqref{eq.kktgrad}--\eqref{eq.kktslack}, then we term $\bX^*$ a \emph{KKT
point} of problem \eqref{eq:symnmf}.

A \emph{stationary point} for problem \eqref{eq:symnmf} is a point $\bX^*$
that satisfies the following optimality condition
\cite[Proposition~2.1.2]{bertsekas99}:
\begin{equation}\label{eq.stationarypoints}
\bigg\langle\big(\bX^*(\bX^*)^{\T}-\frac{\bZ^{\T}+\bZ}{2}\big)\bX^*,\bX-\bX^*\bigg\rangle\ge0,\quad\forall\;\bX\ge0.
\end{equation}

It can be checked that when $\tau$ in \eqref{eq:symnmf:3} is sufficiently
large (larger than a threshold dependent on $\bZ$), then problem
\eqref{eq:symnmf:3} is {\it equivalent} to problem \eqref{eq:symnmf}, in the
sense that the KKT points $\bX^*$ of the two problems are identical. Also,
there is a one-to-one correspondence between the KKT points and stationary
points of the SymNMF problem, although in general such one-to-one
correspondence may not hold. To be more precise, we have:
\begin{lemma}\label{le.kkteqstationary}
For problem \eqref{eq:symnmf}, a point $\bX^*$, is a KKT point,
which means there exists some $\bOmega^*$ such that $(\bX^*, \bOmega^*)$
satisfies \eqref{eq.kktgrad}--\eqref{eq.kktslack}, if and only if $\bX^*$ is a
stationary point, which means it satisfies \eqref{eq.stationarypoints}.
\end{lemma}
\begin{IEEEproof}
See Section~\ref{le.kktandstationary}
\end{IEEEproof}

\begin{lemma}\label{le.bdtau} Suppose $\tau>\theta_k,\forall k$
where
\begin{equation}\label{eq.thetak}
\theta_k\bydef\frac{\bZ_{k,k}+\frac{1}{2}\sqrt{\sum^N_{i=1}(\bZ_{i,k}+\bZ_{k,i})^2}}{2},
\end{equation}
then the KKT points of problem \eqref{eq:symnmf} and the KKT points of
problem \eqref{eq:symnmf:3} have a one-to-one correspondence.
\end{lemma}

\begin{IEEEproof}
See Section~\ref{le.globalopttao}.
\end{IEEEproof}
We remark that the previous work \cite{vagi16} has made the observation that
solving SymNMF with the additional constraints
$\|\bX_i\|_2\le\sqrt{2\|\bZ\|_F},\forall i$ will not result in any loss of the
{\it global} optimality. \leref{le.bdtau} provides a stronger result, that all
{\it KKT} points of SymNMF are preserved within a {\it smaller} bounded
feasible set $\mathcal{Y}\bydef\{\bY\mid
\bY_i\ge0,\|\bY_i\|^2_2\le\tau,\forall i\}$ (note, that $\tau \ll 2\|\bZ\|_F$
in general).

The proposed NS-SymNMF algorithm alternates between the primal updates of
variables $\bX$ and $\bY$, and the dual update for $\bLambda$. Below we
present its detailed steps (superscript $t$ is used to denote the iteration
number).
\begin{align}
\nonumber
\bY^{(t+1)}=&\arg\min_{\bY\ge 0,\|\bY_i\|^2_2\le \tau, \forall
i}\frac{1}{2}\|\bX^{(t)}\bY^{\T}-\bZ\|^2_F
\\
+&\frac{\rho}{2}\|\bY-\bX^{(t)}+\bLambda^{(t)}/\rho\|^2_F+\frac{\beta^{(t)}}{2}\|\bY-\bY^{(t)}\|^2_F,\label{eq.pupdatey}
\\\nonumber
\bX^{(t+1)}=&\arg\min_{\bX}\frac{1}{2}\|\bX(\bY^{(t+1)})^{\T}-\bZ\|^2_F
\\
+&\frac{\rho}{2}\|\bX-\bLambda^{(t)}/\rho-\bY^{(t+1)}\|^2_F,\label{eq.pupdatex}
\\
\bLambda^{(t+1)}=&\bLambda^{(t)}+\rho(\bY^{(t+1)}-\bX^{(t+1)}), \label{eq.pupdatelam}
\\
\beta^{(t+1)}=&\frac{6}{\rho}\|\bX^{(t+1)}(\bY^{(t+1)})^{\T}-\bZ\|^2_F \label{eq.pupdatelab}.
\end{align}
We remark that this algorithm is very close in form to the standard ADMM
method applied to problem \eqref{eq:symnmf:3} (which lacks convergence
guarantees). The key difference is the use of the proximal term
$\|\bY-\bY^{(t)}\|_F^2$ multiplied by an {\it iteration dependent} penalty
parameter $\beta^{(t)}\ge0$, whose value is proportional to the size of the
objective value. Intuitively, if the algorithm converges to a solution with a
small objective value, then parameter $\beta^{(t)}$ vanishes in the limit.
Introducing such proximal term is one of the main novelty of the algorithm,
and it is crucial in guaranteeing the convergence of NS-SymNMF.

\section{Convergence Analysis}

In this section we provide convergence analysis of NS-SymNMF for a general
SymNMF problem. We do not require $\bZ$ to be symmetric,
positive-semidefinite, or to have positive entries. We assume $K$ can be any
integer in $[1, \; N]$.

\subsection{Convergence and Convergence Rate}

Below we present our first main result, which asserts that when the penalty
parameter $\rho$ is sufficiently large, the NS-SymNMF algorithm converges
globally to the set of KKT points of problem \eqref{eq:symnmf}.

\begin{theorem}\label{th.convergence}
Suppose the following is satisfied
\begin{equation}\label{eq.reqrho}
\rho> 6N\tau.
\end{equation}
Then the following statements are true for NS-SymNMF:
\begin{enumerate}

\item The equality constraint is satisfied in the limit, i.e.,
$$\lim_{t\to\infty} \|\bX^{(t)}-\bY^{(t)}\|^2_F\to 0.$$

\item The sequence $\{\bX^{(t)}, \bY^{(t)} \bLambda^{(t)}\}$ generated by the
algorithm is bounded. And every limit point of the sequence is a KKT point of
problem \eqref{eq:symnmf}.
\end{enumerate}
\end{theorem}

An equivalent statement on the convergence is that the sequence $\{\bX^{(t)},
\bY^{(t)} \bLambda^{(t)}\}$ converges to the set of KKT points of problem
\eqref{eq:symnmf}; cf.~footnote~\ref{f.con} on Page~\pageref{f.con}.

\begin{IEEEproof}
See Section~\ref{sec:conproof}.
\end{IEEEproof}

Our second result characterizes the convergence rate of the algorithm. To this
end, we construct a function that measures the optimality of the iterates
$\{\bX^{(t)}, \bY^{(t)}, \bLambda^{(t)}\}$. Define the {\it proximal gradient}
of the augmented Lagrangian function as
\begin{equation}\notag
\widetilde{\nabla}\mathcal{L}(\bX,\bY,\bLambda)\bydef\left[\begin{array}{l}\bY^{\T}-\textsf{proj}_{\mathcal{Y}}[\bY^{\T}-\nabla_{\bY}(\mathcal{L}(\bY,\bX,\bLambda)]\\ \nabla_{\bX}\mathcal{L}(\bX,\bY,\bLambda)\end{array}\right]
\end{equation}
where
\begin{equation}\label{eq.projection}
\textrm{proj}_{\mathcal{Y}}(\bW)\bydef\arg\min_{\bY\ge0,\|\bY_i\|^2_2\le\tau,\forall
i}\|\bW-\bY\|^2_F
\end{equation}
i.e., it is the projection operator that projects a given matrix $\bW$ onto
the feasible set of $\bY$. Here we propose to use the following quantity to
measure the progress of the algorithm
\begin{multline}\label{eq.optgap}
\mathcal{P}(\bX^{(t)},\bY^{(t)},\bLambda^{(t)})\bydef\|\widetilde{\nabla}\mathcal{L}(\bX^{(t)},\bY^{(t)},\bLambda^{(t)})\|^2_F
\onlytwo{\hspace*{10pt}\\} +\|\bX^{(t)}-\bY^{(t)}\|^2_F.
\end{multline}
It can be verified that {if
$\lim_{t\to\infty}\mathcal{P}(\bX^{(t)},\bY^{(t)},\bLambda^{(t)})=0$, then a
KKT point of problem \eqref{eq:symnmf} is obtained.}

Below we show that the function
$\mathcal{P}(\bX^{(t)},\bY^{(t)},\bLambda^{(t)})$ goes to zero in a sublinear
manner.
\begin{theorem}\label{th.convergencerate}
For a given small constant $\epsilon$, let $T(\epsilon)$ denote the iteration
index satisfying the following inequality
\begin{equation}
T(\epsilon)\bydef\min\{t\mid \mathcal{P}(\bX^{(t)},\bY^{(t)},\bLambda^{(t)})\le\epsilon,t\ge0\}.
\end{equation}
Then there exists some constant $C>0$ such that
\begin{equation}
\epsilon\le\frac{C\mathcal{L}(\bX^{(1)},\bY^{(1)},\bLambda^{(1)})}{T(\epsilon)}.
\end{equation}
\end{theorem}
\begin{IEEEproof}
See Section~\ref{sec:conrateproof}.
\end{IEEEproof}
The result indicates that it takes $\mathcal{O}(1/\epsilon)$ iterations for
$\mathcal{P}(\bX^{(t)},\bY^{(t)},\bLambda^{(t)})$ to be less than $\epsilon$.
It follows that NS-SymNMF converges sublinearly.

\subsection{Sufficient Global and Local Optimality Conditions}

Since problem \eqref{eq:symnmf} is not convex, the KKT points obtained by
NS-SymNMF could be different from the global optimal solutions. Therefore it
is important to characterize the conditions under which these two different
types of solutions coincide. Below we provide an easily checkable sufficient
condition to ensure that a KKT point $\bX^*$ is also a globally optimal
solution for problem \eqref{eq:symnmf}.

\begin{theorem}\label{th.globalopt}
Suppose that $\bX^*$ is a KKT point of problem \eqref{eq:symnmf}. Then,
$\bX^*$ is also a global optimal point if the following is satisfied
\begin{equation}\label{eq.reqs}
\mathbf{S}\bydef\bX^*(\bX^*)^{\T}-\frac{\bZ^{\T}+\bZ}{2}\succeq0.
\end{equation}
\end{theorem}
\begin{IEEEproof}
See Section~\ref{sec.globaloptproof}.
\end{IEEEproof}

It is important to note that condition \eqref{eq.reqs} is only a sufficient
condition and hence may be difficult to satisfy in practice. In this section
we provide a milder condition which ensures that a KKT point is {\it locally
optimal}. This type of result is also very useful in practice since it can
help identify spurious saddle points such as the point $\bX^*= \mathbf{0}$ in
the case where $\bZ^{\T}+\bZ$ is not negative semidefinite.

We have the following characterization of the local optimal solution of the
SymNMF problem. \vspace{-0.2cm}
\begin{theorem}\label{th.localopt}
Suppose that $\bX^*$ is a KKT point of problem \eqref{eq:symnmf}. Define a
block matrix $\mathbf{T}\in\mathbb{R}^{KN\times KN}$ whose $(m,n)$th block is
a matrix of size $N \times N$ as follows
\begin{equation}
\mathbf{T}_{m,n}\bydef\left((\bX'^*_m)^{\T}\bX'^*_n-\delta\|\bX'^*_n\|^2_2\right)\bI+\bX'^*_n(\bX'^*_m)^{\T}+ \delta_{m,n}\mathbf{S},
\end{equation}
where $\mathbf{S}$ is defined in \eqref{eq.reqs}, $\delta_{m,n}$ is the
Kronecker delta function, and $\bX'^*_m$ denotes the $m$th column of $\bX^*$.
If there exists some $\delta>0$ such that $\mathbf{T}\succ 0$, then $\bX^*$
is a strict local minimum solution of problem \eqref{eq:symnmf}, meaning that there
exists some $\epsilon>0$ small enough such that for all $\bX\ge0$ satisfying $
\|\bX-\bX^*\|_F\le\epsilon$, we have
\begin{equation}
f(\bX)\ge f(\bX^*)+\frac{\gamma}{2}\|\bX-\bX^*\|^2_F.
\end{equation}
Here the constant $\gamma$ is given by
\begin{equation}\label{eq:gamma}
\gamma=-\left(\frac{2K^2}{\delta}+K(K-2)\right)\epsilon^2+2\lambda_{\min}(\mathbf{T})>0
\end{equation}
where $\lambda_{\min}(\mathbf{T})>0$ is the smallest eigenvalue of
$\mathbf{T}$.
\end{theorem}
\begin{IEEEproof}
See Section~\ref{sec.localoptproof}.
\end{IEEEproof}

In the special case of $K=1$, the sufficient condition set forth in Theorem
\ref{th.localopt} can be significantly simplified.
\begin{corollary}\label{co.localopt}
Suppose that $\bx^*$ is the KKT point of problem \eqref{eq:symnmf} when $K=1$.
If there exists some $\delta>0$ such that
\begin{equation}\label{eq.localoptcond}
\mathbf{T}_1\bydef(1-\delta)\|\bx^*\|^2_2\bI+2\bx^*(\bx^*)^{\T}-\frac{\bZ^{\T}+\bZ}{2}\succ 0,
\end{equation}
then $\bx^*$ is a strict local minimum point of problem \eqref{eq:symnmf}.
\end{corollary}

\begin{IEEEproof}
See Section~\ref{sec.localoptproofkequal1}.
\end{IEEEproof}

We comment that the condition given in Theorem \ref{th.localopt} is much
milder than that in Theorem \ref{th.globalopt}. Further such condition is also
very easy to check as it only involves finding the smallest eigenvalue of a
$KN\times KN$ matrix for a given $\delta$ \footnote{To find such smallest
eigenvalue, we can find the largest eigenvalue of $\eta\bI-\mathcal{T}$, using
algorithms such as the power method
\cite{luwa15}, where $\eta$ is sufficient large based on $\tau$ and
$\|\bZ\|_F$.}. In our numerical results (to be presented shortly), we set
 a series of consecutive $\delta$ when performing the test. We have
observed that the solutions generated by NS-SymNMF satisfy the condition
provided in Theorem \ref{th.localopt} with high probability.

\section{Implementation} \label{sec:majhead}

In this section we discuss the implementation of the proposed algorithm.

\subsection{The $\bX$-Subproblem}

The subproblem for updating $\bX^{(t+1)}$ in \eqref{eq.pupdatex} is equivalent
to the following problem
\begin{equation}
\min_{\bX}\|\bZ^{(t+1)}_{\bX}-\bX\bA^{(t+1)}_{\bX}\|^2_F \label{eq.bigupdatex}
\end{equation}
where
\begin{align}\label{eq.combz}
&\bZ^{(t+1)}_{\bX}\bydef\bZ\bY^{(t+1)}+\bLambda^{(t)}+\rho\bY^{(t+1)}\\
&\bA^{(t+1)}_{\bX}\bydef(\bY^{(t+1)})^{\T}\bY^{(t+1)}+\rho\bI \succ 0\nonumber
\end{align}
are two fixed matrices. Clearly problem \eqref{eq.bigupdatex} is just a
least squares problem and can be solved in closed-form. The solution is given
by
\begin{equation}\label{eq.updatexclosed}
\bX^{(t+1)}=\bZ^{(t+1)}_{\bX}(\bA^{(t+1)}_{\bX})^{-1}.
\end{equation}
We remark that the $\bA^{(t+1)}_{\bX}$ is a $K\times K$ matrix, where $K$ is
usually small (e.g., the number of clusters for graph clustering
applications). As a result, $\bX^{(t+1)}$ in \eqref{eq.updatexclosed} can be
obtained by solving a small system of linear equations and hence
computationally cheap.

\subsection{The $\bY$-Subproblem} \label{ssec:subhead}

The $\bY$-subproblem \eqref{eq.pupdatey} can be decomposed into $N$ separable
constrained least squares problems, each of which can be solved independently,
and hence can be implemented in parallel. We may use the conventional gradient
projection (GP) for solving each subproblem, using iterations
\begin{equation}\label{eq.updatey}
\bY_i^{(r+1)}=\textsf{proj}_{\mathcal{Y}}(\bY_i^{(r)}-\alpha(\bA^{(t)}_{\bY}\bY_i^{(r)}-\bZ^{(t)}_{\bY,i}))
\end{equation}
where
\begin{align}\label{eq.comzy}
\bZ^{(t)}_{\bY}&\bydef(\bX^{(t)})^{\T}\bZ+\rho(\bX^{(t)})^{\T}-(\bLambda^{(t)})^{\T}+\beta^{(t)}(\bY^{(t)})^{\T},
\\
\bA^{(t)}_{\bY}&\bydef(\bX^{(t)})^{\T}\bX^{(t)}+(\rho+\beta^{(t)})\bI\succ 0,
\end{align}
$\bZ_{\bY,i}$ denotes the $i$th column of matrix $\bZ_{\bY}$, $\alpha$ is the
step size, which is chosen either as a constant
$1/\lambda_{\max}(\bA^{(t)}_{\bY})$, or by using some line search procedure
\cite{bertsekas99}; {$r$ denotes the iteration of the inner loop;} for a given
vector $\bw$ , $\textsf{proj}_{\mathcal{Y}}(\bw)$ denotes the projection of it
to the feasible set of $\bY_i$, which can be evaluated in closed-form
\cite[pp.~80]{fapa03} as follows
\begin{align}
\bw^{+}&=\textsf{proj}_+(\bw)\bydef\max\{\bw,\boldsymbol{0}_{K\times 1}\},  \label{eq.proj+}
\\\nonumber
 \bY_i&=\textsf{proj}_{\|\bw^{+}\|^2_2\le\tau}(\bw^+)
 \\
 &\bydef\sqrt{\tau}\bw^+/\max\{\sqrt{\tau},\|\bw^+\|_2\}.\label{eq.projb}
\end{align}
Other algorithms such as accelerated version of the gradient projection
\cite{Beck:2009:FIS:1658360.1658364} can also be used to solve the
$\bY$-subproblem. It is also worth noting that when $\bZ$ is sparse, the
complexity of computing $\bZ\bY^{(t+1)}$ in \eqref{eq.combz} and
$(\bX^{(t)})^{\T}\bZ$ in \eqref{eq.comzy} is only proportional to the number
of nonzero entries of $\bA$.

\section{Numerical Results} \label{sec:print}

\begin{figure*}[htb]
\centering \subfigure[$N=500$, $K=60$.]{
\label{fig:globalopt}
\includegraphics[width=.4\linewidth]{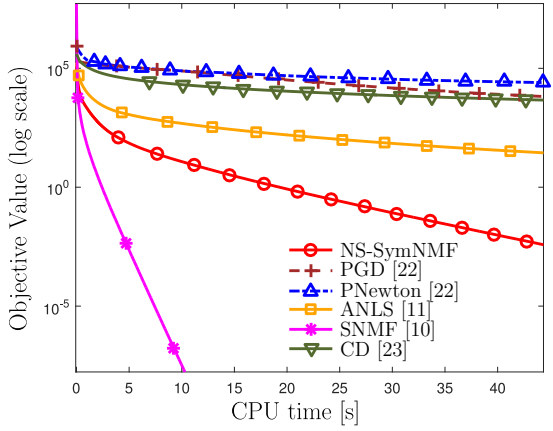}}
\hspace{0.4in} \subfigure[$N=500$, $K=60$, and $\bZ$ is a full rank matrix.]{
\label{fig:localoptk4}
\includegraphics[width=.4\linewidth]{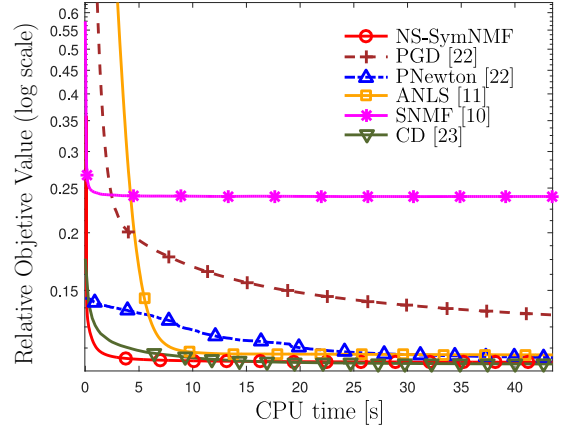}}
\caption{Data Set I: the convergence behaviors of different SymNMF solvers;
each point in the figures is an average of 20 independent MC trials.}
\label{fig:dataseti}
\end{figure*}

\begin{figure*}[htb]
\centering \subfigure[Objective Value]{
\label{fig:syndataobj}
\includegraphics[width=.4\linewidth]{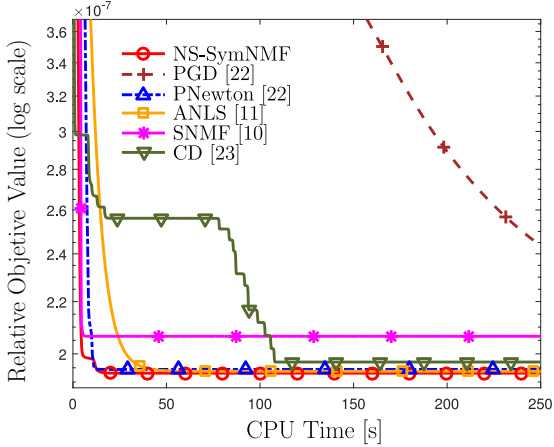}}
\hspace{0.4in} \subfigure[Optimality Gap]{
\label{fig:syndatagap}
\includegraphics[width=.4\linewidth]{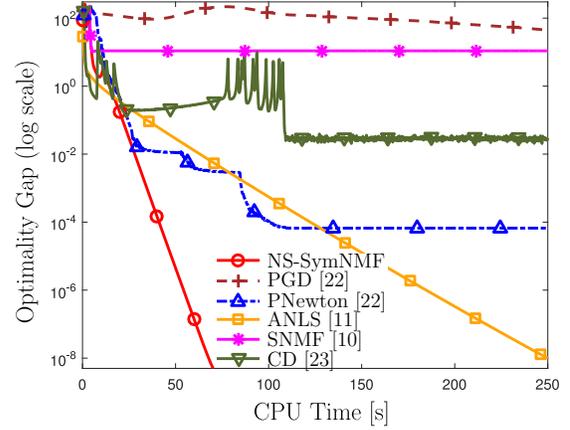}}
\caption{Data Set II: the convergence behaviors of different SymNMF solvers;
each point in the figures is an average of 20 independent MC trials; $N=2000$,
$K=4$}
\label{fig:syndata}
\end{figure*}

In this section, we compare the proposed algorithm with a few existing SymNMF
solvers on both synthetic and real data sets. We run each algorithm with 20
random initializations (except for SNMF, which does not require external
initialization). The entries of the initialized $\bX$ (or $\bY$) follow an
\emph{i.i.d.} uniform distribution in the range $[0,\tau]$. All algorithms are
started with the same initial point each time, and all tests are performed
using Matlab on a computer with Intel Core i5-5300U CPU running at 2.30GHz
with 8GB RAM. Since the compared algorithms have different computational
complexity, we use the objective values versus CPU time for fair comparison.
We next describe different SymNMF solvers that are compared in our work.

\noindent{\bf Algorithms Comparison.} In our numerical simulations, we compare
the following algorithms.

\paragraph{Projected Gradient Descent (PGD) and Projected Newton method
(PNewton) \cite{kudi12,kuyu15}} The PGD and PNewton directly use the gradient
of the objective function. The key difference between them is that PGD adopts
the identity matrix as a scaling matrix while PNewton exploits reduced Hessian
for accelerating the convergence rate. The PGD algorithm converges slowly if
the step size is not well selected, while the PNewton algorithm has high
per-iteration complexity compared with ANLS and NS-SymNMF, due to the
requirement of computing the Hessian matrix. Note that to the best of our
knowledge, neither PGD nor PNewton possesses convergence or rate of
convergence guarantees.

\paragraph{Alternating Nonnegative Least Square (ANLS) \cite{kuyu15}} The ANLS
method is a very competitive SymNMF solver, which can be implemented in
parallel easily. ANLS reformulates SymNMF as
\begin{equation}\notag
\min_{\bX,\bY\ge0}g(\bX,\bY) = \|\bX\bY^{\T}-\bZ\|^2_F+\nu\|\bX-\bY\|^2_F
\end{equation}
where $\nu>0$ is the regularization parameter. One of shortcomings is that
there is no theoretical guarantee that the ANLS method can converge to the set
of KKT points of problem \eqref{eq:symnmf} or even producing two symmetric
factors, although a penalty term for the difference between the factors ($\bX$
and $\bY$) is included in the objective.

\paragraph{Symmetric Nonnegative Matrix Factorization (SNMF) \cite{kesi14}}
The SNMF algorithm transforms the original problem to another one under the
assumption that $\bZ$ can be exactly decomposed by $\bX\bX^{\T}$. Although
SNMF often converges quickly in practice, there has been no theoretical
analysis under the general case where $\bZ$ cannot be exactly decomposed.

\paragraph{Coordinate Descent (CD) \cite{vagi16}} The CD method updates each
entry of $\bX$ in a cyclic way. For updating each entry, we only need to find
the roots of a fourth-order univariate function. However, CD may not converge
to the set of KKT points of SymNMF. Instead, there is an additional condition
given in \cite{vagi16} for checking whether the generated sequence converges
to a unique limit point. A heuristic method for checking the condition is
additionally provided, which requires, e.g., plotting the norm between the
different iterates.

\paragraph{The Proposed NS-SymNMF} The update rule of NS-SymNMF is similar to
that of ANLS. The difference between them is that NS-SymNMF uses one
additional block for dual variables and ANLS adds a penalty term. The dual
update involved in NS-SymNMF benefits the convergence of the algorithm to KKT
points of SymNMF.

We remark that in the implementation of NS-SymNMF we let $\tau=\max_k \theta_k
$ (cf.~\eqref{eq.thetak}) and the maximum number of iterations of GP be $40$.
Also, we gradually increase the value of $\rho$ from an initial value to meet
condition \eqref{eq.reqrho} for accelerating the convergence
rate \cite{raho14}. Here, the choice of $\rho$ follows $\rho^{(t+1)} =
\min\{\rho^{(t)}/(1-\epsilon/\rho^{(t)}),6.1N\tau\}$ where $\epsilon=10^{-3}$
as suggested in \cite{scfa14}. We choose $\rho^{(1)}=\bar{\tau}$ for the case
that $\bZ$ can be exactly decomposed and $\sqrt{N}\bar{\tau}$ for the rest of
cases, where $\bar{\tau}$ is the mean of $\theta_k,\forall k$. The similar
strategy is also applied for updating $\beta^{(t)}$. We choose
$\beta^{(t)}=6\xi^{(t)}\|\bX^{(t)}\bY^{(t)}-\bZ\|^2_F/\rho^{(t)}$ where
$\xi^{(t+1)}=\min\{\xi^{(t)}/(1-\epsilon/\xi^{(t)}),1\}$ and $\xi^{(1)}=0.01$,
and only update $\beta^{(t)}$ once every 100 iterations to save CPU time. To
update $\bY$, we implement the block pivoting method
\cite{kipa11} since such method is faster than the GP method for solving the
nonnegative least squares problem. If $\|\bY^{(t+1)}_i\|^2_2\le\tau$ is not
satisfied, then we switch to GP on $\bY^{(t)}_i$. We also remark
that we set the step size of PGD to $10^{-5}$ for all tested cases, and use
the Matlab codes of PNewton and ANLS from
\url{http://math.ucla.edu/~dakuang/}.

\noindent{\bf Performance on Synthetic Data.} First we describe the two synthetic data
sets that we have used in the first part of the numerical results.

\noindent \underline{Data set I (Random symmetric matrices):} We
randomly generate two types of symmetric matrices, one is of low rank and the
other is of full rank.

For the low rank matrix, we first generate a matrix $\bM$ with dimension
$N\times K$, whose entries follow an \emph{i.i.d.} Gaussian distribution with
zero mean and unit variance. We use $\bM_{i,j}$ to denote the $(i,j)$th entry
of $\bM$. Then generate a new matrix $\widetilde{\bM}$ whose $(i,j)$th entry
is $|\bM_{i,j}|$. Finally, we obtain a positive symmetric
$\bZ=\widetilde{\bM}\widetilde{\bM}^{\T}$ as the given matrix to be
decomposed.

For the full rank matrix, we first randomly generate a $N\times N$ matrix
$\bP$, whose entries follow an \emph{i.i.d.} uniform distribution in the
interval $[0,1]$. Then we compute $\bZ=(\bP+\bP^{\T})/2$.

\noindent\underline{Data set II (Adjacency matrices):} One important
application of SymNMF is graph partitioning, where the adjacency matrix of a
graph is factorized. We randomly generate a graph as follows. First, set the
number of nodes to $N$ and the number of cluster to $4$, and the numbers of
nodes within each cluster to $300,500,800,400$. Second, we randomly generate
data points whose relative distance will be used to construct the adjacency
matrix. Specifically, data points $\{x_i\}\in\mathbb{R}$, $i=1,\ldots,N$, are
generated in one dimension. Within one cluster, data points follow an
\emph{i.i.d.} Gaussian distribution. The means of the random variables in
these 4 clusters are $2,3,6,8$, respectively, and the variance is 0.5 for all
distributions. Construct the similarity matrix $\bA\in\mathbb{R}^{N\times N}$,
whose $(i,j)$th entry is $\bA_{i,j}=\exp(-(x_i-x_j)^2/(2\sigma^2))$ where
$\sigma^2=0.5$.

The convergence behaviors of different SymNMF solvers for the synthetic data
sets are shown in \figref{fig:dataseti} and \figref{fig:syndata}. The results
are averaged over 20 Monte Carlo (MC) trials with independently generated
data. In \figref{fig:globalopt}, the generated $\bZ$ can be exactly decomposed
by SymNMF. It can be observed that NS-SymNMF and SNMF converge to the global
optimal solution quickly, and SNMF is the fastest one among all compared
algorithms. However, the case where the matrix can be exactly factorized is
not common in most practical applications. Hence, we also consider the case
where matrix $\bZ$ cannot be factorized exactly by a $N\times K$ matrix. The
results are shown in \figref{fig:localoptk4} and we use the relative objective
value for comparison, i.e., $\|\bX\bX^{\T}-\bZ\|^2_F/\|\bZ\|^2_F$. We can
observe that NS-SymNMF and CD can achieve a lower objective value than other
methods. It is worth noting that there is a gap between SNMF and others, since
the assumption of SNMF is not satisfied in this case.

We also implement the algorithms on the adjacency matrices (data set II),
where the results are shown in \figref{fig:syndata}. The NS-SymNMF and SNMF
algorithms converge very fast, but it can be observed that there is still a
gap between SNMF and NS-SymNMF as shown in \figref{fig:syndataobj}. We further
show the convergence rates with respective to optimality gap versus CPU time
in \figref{fig:syndatagap}. The optimality gap \eqref{eq.optgap} measures the
closeness between the generated sequence and the true stationary point. To get
rid of the effect of the dimension of $\bZ$, we use
$\|\bX-\textsf{proj}_+[\bX-\nabla_{\bX}(f(\bX))]\|_{\infty}$ as the optimality
gap. It is interesting to see the ``swamp'' effect \cite{swamp}, where the
objective value generated by the CD algorithm remains almost constant during
the time period from around 25s to 75s although actually the corresponding
iterates do not converge, and then the objective value starts decreasing
again.

\begin{figure}[htb]
\centering \includegraphics[width=\onlytwo{.8}\onlyone{.5}\linewidth]{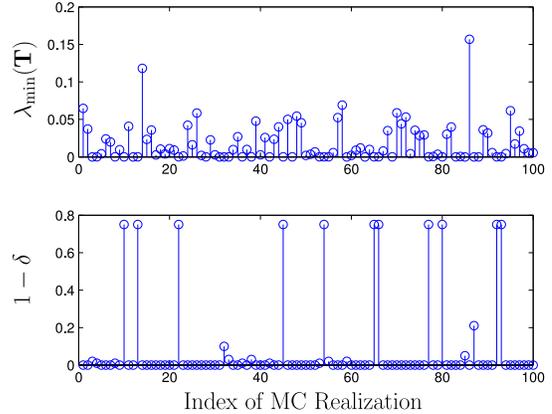}
\caption{Checking local optimality condition, where $N=500$.}
\label{fig:localkkt}
\end{figure}

\noindent{\bf Checking Global/Local Optimality.} After the NS-SymNMF algorithm
has converged, the local/global optimality can be checked according to
\thref{th.globalopt} and \thref{th.localopt}. To find an appropriate $\delta$
that satisfying the condition where $\lambda_{\min}(\mathbf{T})>0$, we initialize
$\delta$ as 1 and decrease it by $0.01$ each time and check the minimum
eigenvalue of $\mathbf{T}$. Here, we use data set II with the fixed ratio of
the number of nodes within each cluster (i.e., $3:5:8:4$) and test on the
different total numbers of nodes. The simulation results are shown in
\figtab{tab:opt} with 100 MC trials, where the average value of
$\lambda_{\min}(\mathbf{T})$ and $\delta$ are given. Further, the percentage
of being able to find a valid $\delta>0$ that ensures $\lambda_{\min}(\mathbf{T})>0$
is listed as the last column. We note that there always existed a $\delta$
such that $\mathbf{T}$ is positive definite in all cases that we tested. This
indicates that (with high probability) the proposed algorithm converges to a
locally optimal solution. In \figref{fig:localkkt}, we provide the values of
$\delta$ that make the corresponding $\lambda_{\min}(\mathbf{T})>0$ at each
realization.

 We also remark that in practice we stop the
algorithm in finite steps, so only an approximate KKT point will be obtained,
and the degree of such approximation can be measured by the optimality gap
defined in \eqref{eq.optgap}.
\begin{table}[htp]
\small \centering
\caption{Local Optimality}\label{tab:opt}
\begin{tabular}{r|c|c|c}
\hline \multicolumn{1}{c}{ $N$ } & \multicolumn{1}{|c|}{
$\lambda_{\min}(\mathbf{T})$} & \multicolumn{1}{c|}{ $\delta$} &
\multicolumn{1}{c}{ Local Optimality (true)}\\ \hline
 50 & $2.71\times 10^{-4}$  & 0.42  & 100\%  \\
\hline
 100 & $4.16\times 10^{-4}$ & 0.37  & 100\% \\
\hline
 500 & $1.8\times 10^{-2}$ &  0.91 & 100\% \\
\hline
\end{tabular}
\end{table}

\begin{figure*}[htb]
  \centering
  \subfigure[Mean of the objective values: Reuters data set]{
    \label{fig:corpora_mean}
    \includegraphics[width=.4\linewidth]{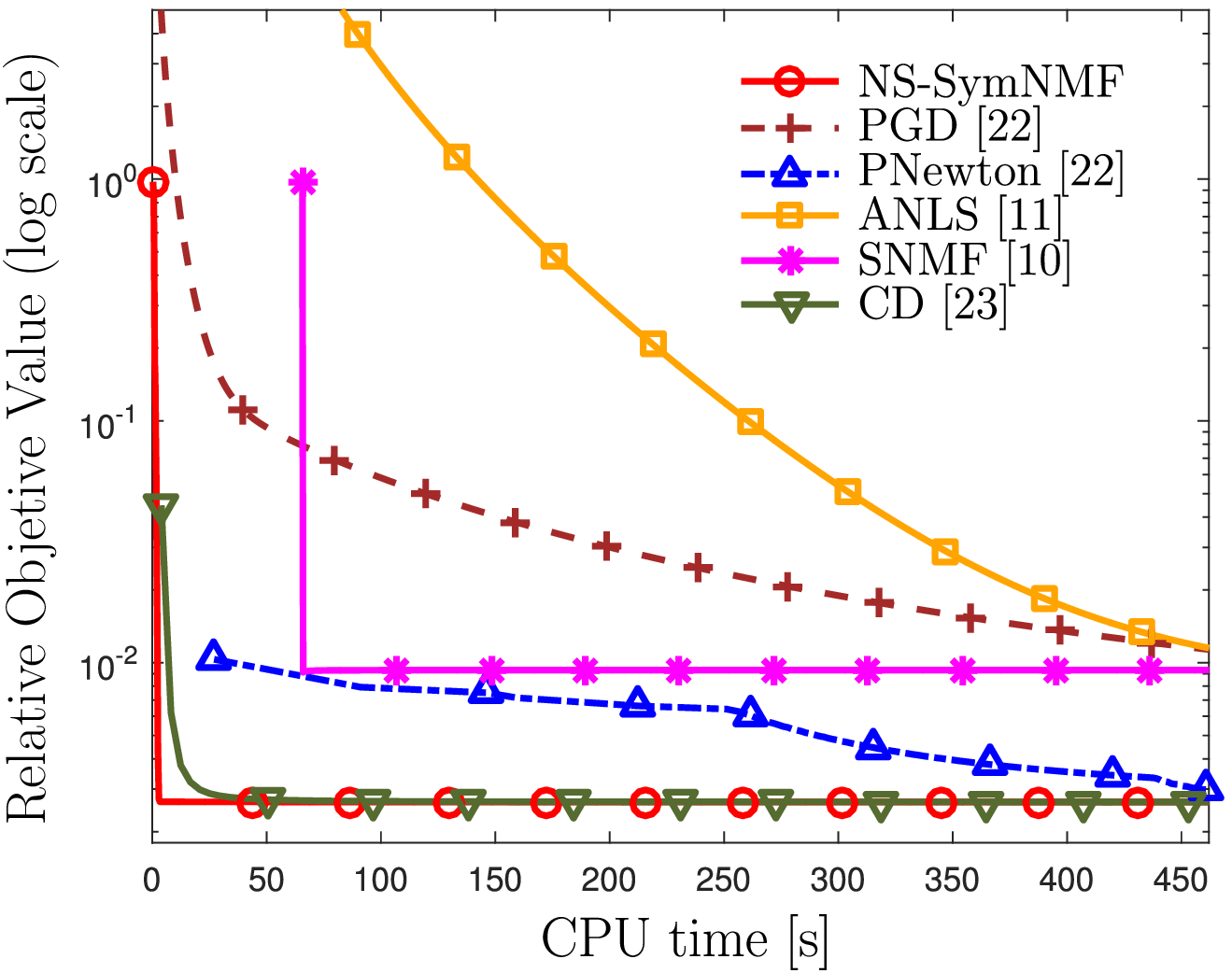}}
  \hspace{0.4in}
  \subfigure[{Mean of the objective values: TDT2 data set}]{
    \label{fig:tdt2}
    \includegraphics[width=.4\linewidth]{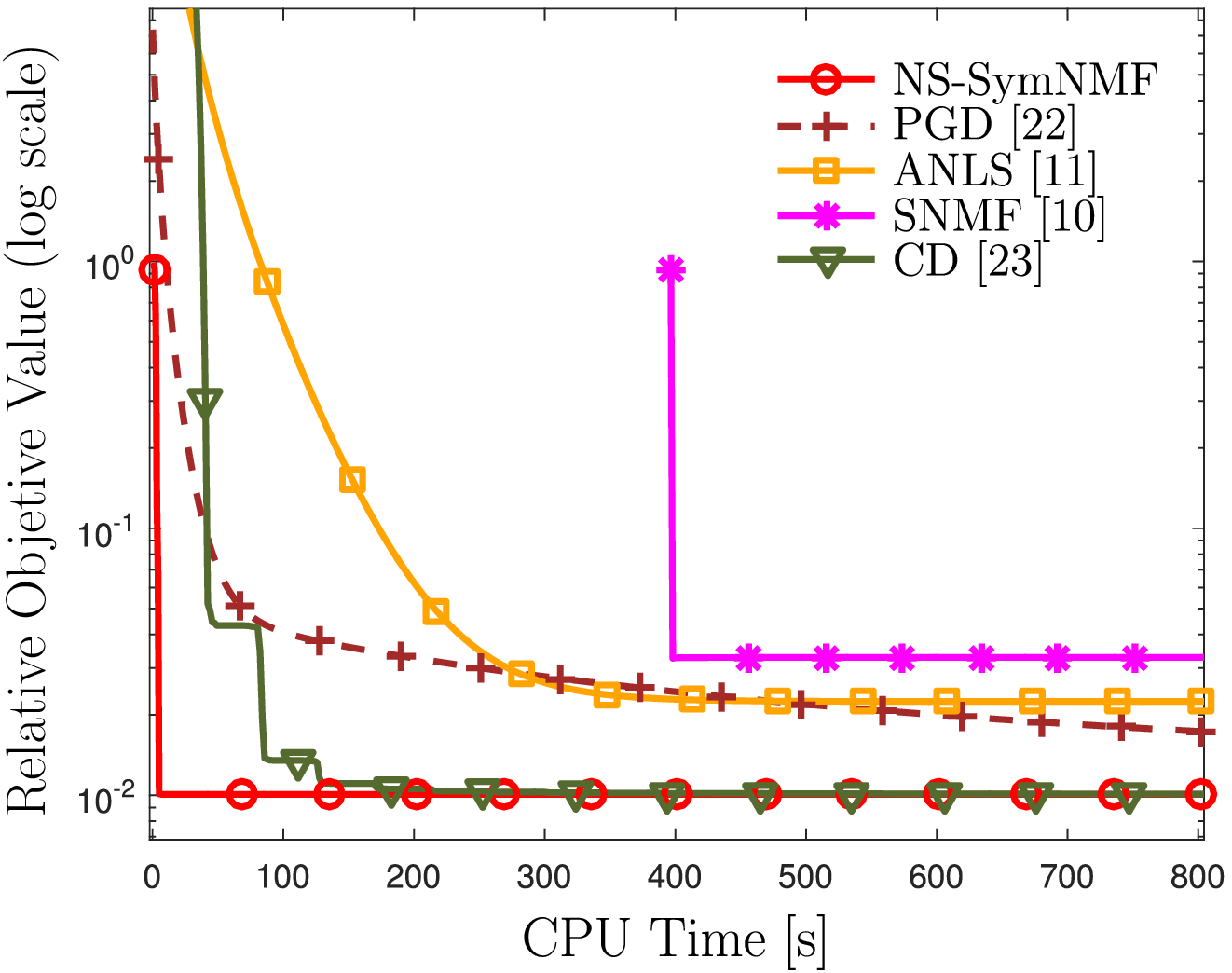}}
  \caption{The convergence behaviors of different SymNMF solvers for the dense similarity matrix; each point in the figures is an average of 20 independent MC trials based on random initializations.}
  \label{fig:densedata}
\end{figure*}

\begin{figure*}[htb]
  \centering
  \subfigure[Mean of the objective values: email-Enron data set]{
    \label{fig:email}
    \includegraphics[width=.4\linewidth]{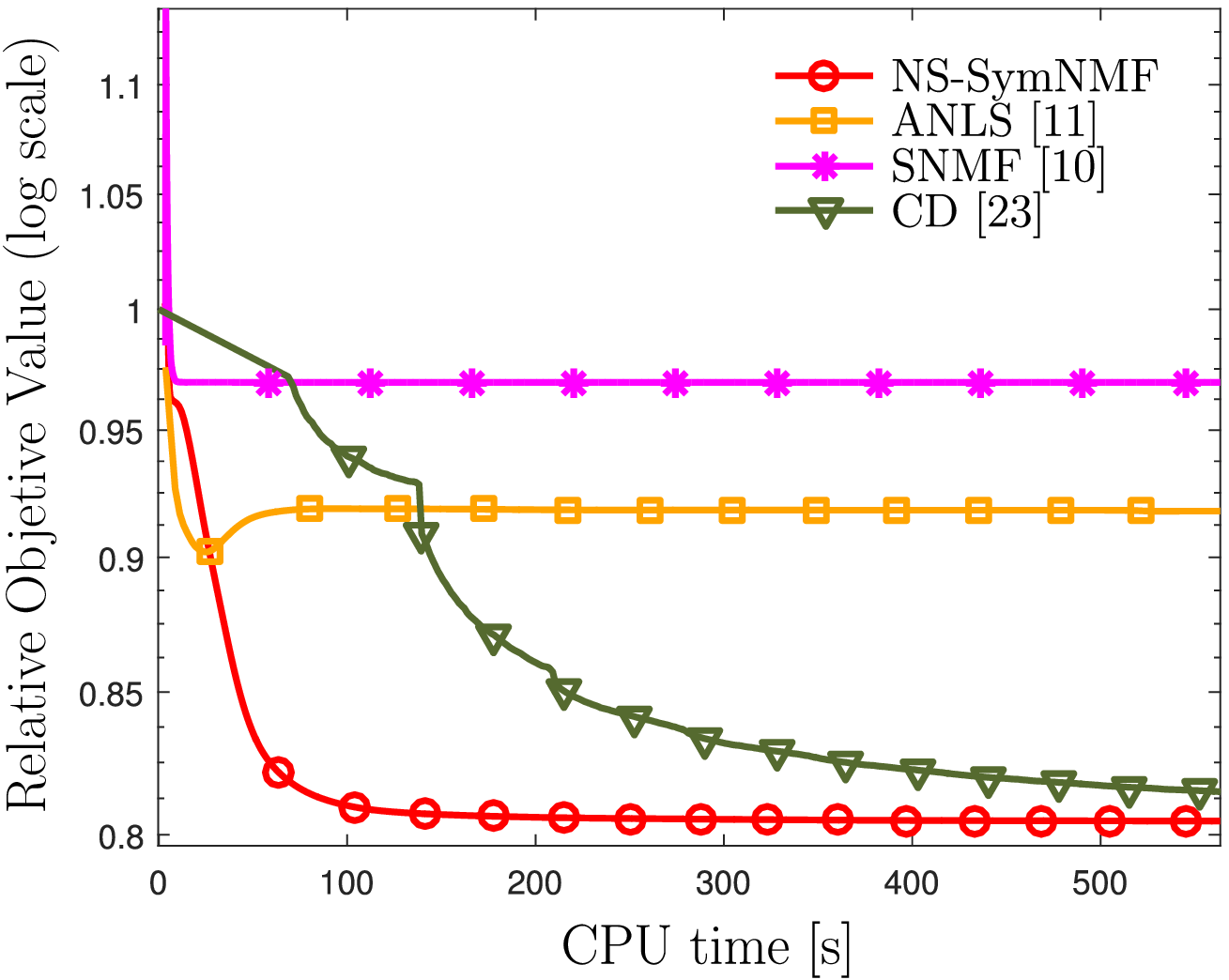}}
  \hspace{0.4in}
  \subfigure[{Mean of the objective values: loc-Brightkite data set}]{
    \label{fig:loc}
    \includegraphics[width=.4\linewidth]{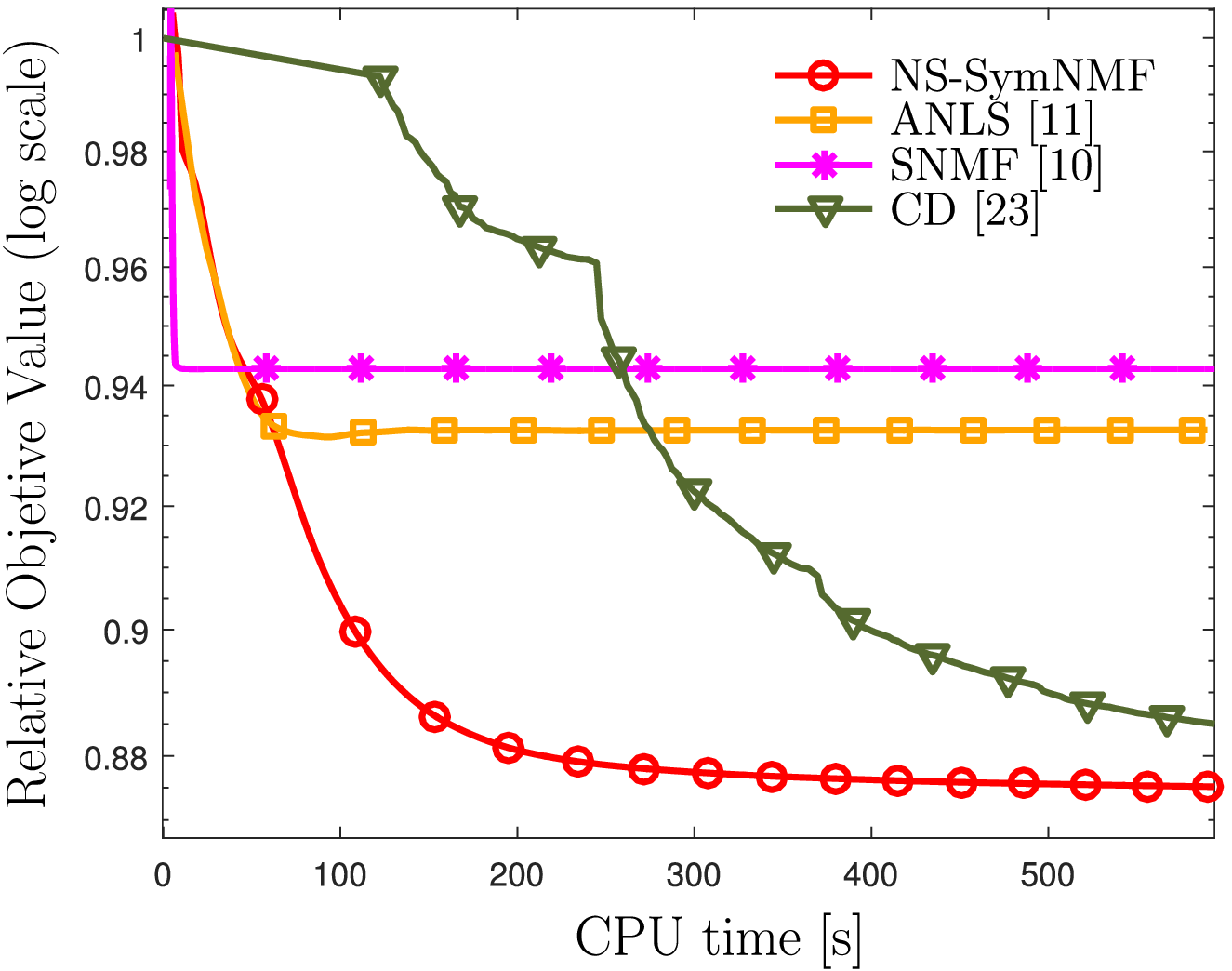}}
  \caption{The convergence behaviors of different SymNMF solvers for the sparse similarity matrix; each point in the figures is an average of 20 independent MC trials based on random initializations.}
  \label{fig:largesparsedata}
\end{figure*}

{\color{blue}
\begin{table*}[htp]
\onlytwo{\footnotesize}\onlyone{\fontsize{2.5mm}{2.5mm}\selectfont }
\centering
\caption{Mean and Standard Deviation of $\|\bX\bX^{\T}-\bZ\|^2_F/\|\bZ\|^2_F$
of the Final Solution of Each Algorithm based on Random
Initializations}\label{tab:dense_meanandstd}
\centering
\begin{tabular}{c||c|c|c|c|c|c|c|c}
\hline Dense Data Sets & $N$　&　$K$ & NS-SymNMF & PGD \cite{kudi12} &
PNewton \cite{kudi12} & ANLS \cite{kuyu15} & SNMF \cite{kesi14} & CD
\cite{vagi16} \\ \hline Reuters
\cite{cahe11} & 4,633 & 25 & {\bf 2.65e{-3}$\pm$3.31e-10} &
1.14e-2$\pm$1.18e-5 & 2.98e-3$\pm$3.71e-6 & 1.16e-2$\pm$1.61e-5 & 9.32e-3 &
2.66e-3$\pm$2.04e-8 \\ \hline TDT2
\cite{cahe11} & 8,939 & 25 & {\bf 1.01e{-2}$\pm$5.35e-9} & 1.74e-2$\pm$7.34e-6
& - & 2.25e-2$\pm$1.25e-6 & 3.29e{-2} & 1.01e-2$\pm$1.21e-6 \\\hline
\end{tabular}
\end{table*}
}

{\color{blue}
\begin{table*}[htp]
\footnotesize
\centering
\caption{Mean and Standard Deviation of $\|\bX\bX^{\T}-\bZ\|^2_F/\|\bZ\|^2_F$
of the Final Solution of Each Algorithm based on Random
Initializations}\label{tab:sparse_meanandstd}
\centering
\begin{tabular}{c||c|c|c|c|c|c|c}
\hline Sparse Data Sets & $N$　&　$K$ & \#nonzero & NS-SymNMF & ANLS
\cite{kuyu15} & SNMF \cite{kesi14} & CD \cite{vagi16} \\ \hline email-Enron
\cite{lela09} &　36,692　& 50 & 367,662 & {\bf 8.05e-1$\pm$4.66e-4} &
9.18e-1$\pm$6.20e-3 & 9.69e-1 &8.13e-1$\pm$1.47e-3 \\ \hline loc-Brightkite
\cite{chmy11} &　58,228　& 50 & 428,156 & {\bf 8.75e-1$\pm$9.52e-4} &
9.33e-1$\pm$1.93e-3 & 9.43e-1 & 8.84e-1$\pm$1.49e-3 \\ \hline
\end{tabular}
\end{table*}
}

\noindent{\bf Performance on Real Data.} We also implement the algorithm on a few real data sets in clustering
applications, which will be described in the next paragraphs.

\subsubsection{Dense Similarity Matrix} we generate the dense
similarity matrices based on the two real data sets: \underline{Reuters-21578
and TDT2
\cite{cahe11}}. We use the 10th subset of the processed Reuters-21578 data
set, which includes $N=4,633$ documents divided into $K=25$ classes. The
number of features is 18,933. Topic detection and tracking 2 (TDT2) corpus
includes two newswires (APW and NYT), two radio programs (VOA and PRI) and two
television programs (CNN and ABC). We use the 10th subset of the processed
TDT2 data set with $K=25$ classes which includes $N=8,939$ documents and each of them
has 36,771 features. We comment that the 10th TDT2 subset is the largest among
the all TDT2 and Reuters subsets. Any other subset can be used equally well.
The similarity matrix is constructed by the Gaussian function where the
difference between two documents is measured by all features using the
Euclidean distance
\cite{cahe11}.

The means and standard deviations of the objective values of the final
solutions are shown in \figtab{tab:dense_meanandstd}. Convergence results of
the algorithms are shown in \figref{fig:densedata}. For the Reuters and TDT2
datasets, before SNMF completes the eigenvalue decomposition for the first
iteration, CD and NS-SymNMF have already obtained low objective values. Also,
since calculating Hessian in PNewton is time consuming, the result of PNewton
is out of range in \figref{fig:tdt2}.

\subsubsection{Sparse Similarity Matrix} we also generate
multiple convergence curves for each algorithm with random initializations
based on some sparse real data sets.

\noindent \underline{Email-Enron network data set \cite{lela09}:} Enron email
corpus includes around half million emails. We use the relationships between
two email addresses to construct the similarity matrix for decomposing. If an
address $i$ sent at least one email to address $j$, then we take
$\bA_{i,j}=\bA_{j,i}=1$. Otherwise, we set $\bA_{i,j}=\bA_{j,i}=0$.

\noindent \underline{Brightkite data set \cite{chmy11}:} Brightkite was a
location-based social networking website. Users were able to share their
current locations by checking-in. The friendships of the users were maintained
by Brightkite. The way of constructing the similarity matrix is the same as
the Enron email data set.

The means and standard deviations of the objective values of the final
solutions are shown in \figtab{tab:sparse_meanandstd}. From the simulation
results shown in \figref{fig:largesparsedata}, it can be observed that the
NS-SymNMF algorithm converges faster than CD, while SNMF and ANLS converge to
some points where the relative objective values are higher than the one
obtained by NS-SymNMF.

\section{Conclusions} \label{sec:conclusion}

In this paper, we propose a nonconvex splitting algorithm for solving the
SymNMF problem. We show that the proposed algorithm converges to a KKT point
in a sublinear manner. Further, we provide sufficient conditions to identify
global or local optimal solutions of the SymNMF problem. Numerical experiments
show that the proposed method can converge quickly to local optimal solutions.

In the future, we plan to extend the proposed methods in a way such that the
algorithms can converge to the local or even global optimal solutions of
SymNMF without requiring checking conditions. Also, it is possible to apply the nonconvex splitting method to more general
matrix factorization problems, such as the quadratic nonnegative matrix
factorization problem \cite{yaaj2012}.

\section{Appendix}

\subsection{Proof of \leref{le.kkteqstationary}}\label{le.kktandstationary}

Sufficiency: the stationary points satisfy
\begin{equation}\label{eq.kkteqstation}
\bigg\langle\big(\bX^*(\bX^*)^{\T}-(\bZ^{\T}+\bZ)/2\big)\bX^*,\bX-\bX^*\bigg\rangle\ge0,\quad\forall\;\bX\ge0.
\end{equation}
Let $\bOmega\bydef(\bX^*(\bX^*)^{\T}-(\bZ^{\T}+\bZ)/2)\bX^*/2$. We have
$\langle\bOmega,\bX-\bX^*\rangle\ge0,\forall\bX\ge0$. By setting $\bX$
appropriately as $0\le \bX\le\bX^*$, we have $\bOmega_{i,j}\ge 0,
(i,j)\in\mathcal{S}$ where $\mathcal{S} = \{i,j|\bX^*_{i,j}\neq0\}$. Also, by
setting $\bX$ appropriately as $\bX\ge \bX^*$, we have $\bOmega_{i,j}\ge 0,
(i,j)\notin\mathcal{S}$. Combining the two cases, we conclude that
$\bOmega\ge0$.

From \eqref{eq.kkteqstation}, we know that $\langle\bOmega,\bX\rangle\ge
\langle\bOmega,\bX^*\rangle$. Since $\bOmega\ge0$ and $\bX\ge0$, we have
$\langle\bOmega,\bX\rangle\ge0,\forall\bX$, meaning that
$\langle\bOmega,\bX^*\rangle\le0$. Combining with $\bX^*\ge0$ and
$\bOmega\ge0$, we have $\langle\bOmega,\bX^*\rangle\ge0$, which results in
$\langle\bOmega,\bX^*\rangle=0$.

In summary, we have
\begin{subequations}
\begin{align}
&2\left(\bX^*(\bX^*)^{\T}-\frac{\bZ^{\T}+\bZ}{2}\right)\bX^*-\bOmega=0,
\\
&\bOmega\ge 0,
\\
&\bX^*\ge 0,
\\
&\langle\bX^*,\bOmega\rangle=0,
\end{align}
\end{subequations}
which are the KKT conditions of the SymNMF problem.

Necessity: If the point is a KKT point of SymNMF, we have
\begin{equation}\label{eq.omegastar}
\bOmega^*=2\big(\bX^*(\bX^*)^{\T}-\frac{\bZ^{\T+}\bZ}{2}\big)\bX^*.
\end{equation}
Combining with $\langle\bX^*,\bOmega^*\rangle=0$, we know that
\begin{equation}
\langle\bOmega^*,\bX-\bX^*\rangle\ge0,\quad\forall\;\bX\ge0,
\end{equation}
which is the condition of stationary points.

\subsection{Proof of \leref{le.bdtau}}
\label{le.globalopttao}

We prove that if $\tau$ is large enough, then the KKT conditions of
\eqref{eq:symnmf} and \eqref{eq:symnmf:3} are the same.
\begin{IEEEproof}
It is sufficient to show that when $\tau$ is large enough, there can be no KKT
point whose column has size $\tau$, leading to the fact that the constraint
$\|\bX^*_k\|^2\le \tau$ is always inactive.

We check the optimality condition of the SymNMF problem at
$\|\bX^*_k\|^2=\tau_k$, where $\tau_k>0$ is a constant. We can rewrite the
objective function as
\begin{align}
\nonumber
f(\bX)=&\frac{1}{2}\bigg(\sum^N_{i=1,i\ne k}\sum^N_{j=1,j\ne k}(\bX_i\bX^{\T}_j-\bZ_{i,j})^2
\onlytwo{\\\notag &}
+\sum^N_{i=1,i\ne k}(\bX_i\bX^{\T}_k-\bZ_{i,k})^2
\\\nonumber
&+\sum^N_{j=1,j\ne k}(\bX_k\bX^{\T}_j-\bZ_{k,j})^2+(\bX_k\bX^{\T}_k-\bZ_{k,k})^2\bigg).
\end{align}
Note, $\bX_i,\bX_j,\bX_k$ denote rows of matrix $\bX$.

We take the gradient of $f(\bX)$ with respective to $\bX_k$:
\begin{align}
\notag
&\frac{\partial f(\bX)}{\partial \bX_{k,m}}=\sum^N_{i=1,i\ne k}\bX_{i,m}(\bX_i\bX^{\T}_k-\bZ_{i,k})
\\\notag
&\sum^N_{j=1,j\ne k}\bX_{j,m}(\bX_k\bX^{\T}_j-\bZ_{k,j})+2\bX_{k,m}(\bX_k\bX^{\T}_k-\bZ_{k,k})
\\\notag
&=\sum^N_{i=1,i\ne k}\bX_{i,m}(\bX_i\bX^{\T}_k-(\bZ_{i,k}+\bZ_{k,i}))
\\
&\qquad\qquad+2\bX_{k,m}(\bX_k\bX^{\T}_k-\bZ_{k,k})
\label{eq.gradientofX}
\end{align}
where $\bX_{i,m}$ denotes the $m$th entry of the $i$th row of $\bX$.

Assume that $\bX^*_k$ is a KKT point. We have $(\frac{\partial
f(\bX^*_k)}{\partial \bX_k})(\bX_k-\bX^*_k)^{\T}\ge 0,
\forall~\bX_k\in\mathcal{X}$, where $\mathcal{X}=\{\bX_k|\bX_k\ge 0,
\|\bX_k\|^2\le\tau_k\}$, which implies
\begin{align}\label{eq.localoptreq1}
&\frac{\partial f(\bX^*_k)}{\partial \bX_{k,m}}(\bX_{k,m}-\bX^*_{k,m})\ge 0
\notag \\
& 0\le\bX_{k,m}\le\bX^*_{k,m}=\sqrt{\tau_k-\sum^K_{n=1,n\ne m}(\bX^*_{k,n})^2}\quad\forall~m.
\end{align}

Since $\|\bX^*_k\|^2=\tau_k$, there exists an index $m$ such that
$\bX^*_{k,m}>0$. Consider a feasible point $0\le\bX_{k,m}<\bX^*_{k,m}$, where
$m\in\mathcal{S}_m\bydef\{m|\bX^*_{k,m}\ne 0\}$. Thanks to
\eqref{eq.localoptreq1}, we have
\begin{equation}\label{eq.optconditionofX}
\frac{\partial f(\bX^*_{k,m})}{\partial \bX_{k,m}}\le0,\quad 0\le\bX_{k,m}<\bX^*_{k,m}\quad\forall~m\in\mathcal{S}_m.
\end{equation}

Plugging \eqref{eq.gradientofX} into \eqref{eq.optconditionofX} and
multiplying $\bX^*_{k,m}$ on both sides of \eqref{eq.optconditionofX}, we can
obtain
\begin{equation}
\def\lhs{\bX^*_{k,m}\bigg(\sum^N_{i=1,i\ne k}\bX^*_{i,m}
  \big(\bX^*_i(\bX^*_k)^{\T}-\frac{\bZ_{i,k}+\bZ_{k,i}}{2}\big)}
\def\rhs{+\bX^*_{k,m}(\bX^*_k(\bX^*_k)^{\T}-\bZ_{k,k})\bigg)\le 0
  \quad\forall~m\in\mathcal{S}_m.}
\iftoggle{twocolumn}{
  \begin{split}
  & \lhs \\ & \qquad \rhs
  \end{split}
}{
  \lhs \rhs
}
\label{eq.optlocalcheck}
\end{equation}

For the case $m\notin\mathcal{S}_m$, we know that $\bX^*_{k,m}=0$. Summing up
\eqref{eq.optlocalcheck} $\forall m$, and noting that $|\mathcal{S}_m|\ge 1$
we can get
\begin{equation}
\def\lhs{p\bydef\sum^N_{i=1,i\ne k}
  \underbrace{\bX^*_{i}(\bX^*_k)^{\T}\big(\bX^*_i(\bX^*_k)^{\T}-\frac{\bZ_{i,k}+\bZ_{k,i}}{2}\big)}
  _{\bydef\mathcal{M}_{i,k}}}
\def\rhs{+\bX^*_{k}(\bX^*_k)^{\T}(\bX^*_k(\bX^*_k)^{\T}-\bZ_{k,k})\le 0.}
  \begin{split}
  \lhs \onlytwo{\\} \rhs
  \end{split}
\label{eq.localptau}
\end{equation}
In \eqref{eq.localptau}, $\mathcal{M}_{i,k}$ is a quadratic function with
respective to $C_{i,k}$, where $C_{i,k}\bydef\bX^*_{i}(\bX^*_k)^{\T}$, so
the minimum of $\mathcal{M}_{i,k}$ is $-1/4((\bZ_{i,k}+\bZ_{k,i})/2)^2$.
Consequently, the minimum of $\sum^N_{i=1,i\ne k}\mathcal{M}_{i,k}$ is
$-1/4\sum^N_{i=1,i\ne k}((\bZ_{i,k}+\bZ_{k,i})/2)^2$.

In addition, since we have $\|\bX^*_k\|^2=\tau_k$, the lower bound of $p$ is
$p_{\textsf{L}}\bydef-1/4\sum^N_{i=1,i\ne
k}((\bZ_{i,k}+\bZ_{k,i})/2)^2+\tau_k(\tau_k-\bZ_{k,k})$ which is a quadratic
function in terms of $\tau_k$. Therefore, if
\begin{equation}
\tau_k>\theta_k\bydef\frac{\bZ_{k,k}+\frac{1}{2}\sqrt{\sum^N_{i=1}(\bZ_{i,k}+\bZ_{k,i})^2}}{2},
\end{equation}
then $p\ge p_{\textsf{L}}>0$, which contradicts the optimality condition
\eqref{eq.optlocalcheck}. It can be concluded that whenever $\tau_k$ is large
enough, at any KKT point no column will have size equal to $\tau_k$.
Furthermore, it can be easily checked that $\tau>\max_k \theta_k$ is a
sufficient condition. The proof is complete.
\end{IEEEproof}

\subsection{Convergence Proof of the Proposed Algorithm} \label{sec:conproof}

In this section, we prove \thref{th.convergence}. The analysis consists of a
series of lemmas.
\begin{lemma}\label{le.bdofdual}
Consider using the update rules \eqref{eq.pupdatey} -- \eqref{eq.pupdatelam}
to solve problem \eqref{eq:symnmf}. Then we have
\begin{align}
\nonumber
\|\bLambda^{(t+1)}-&\bLambda^{(t)}\|^2_F\le3N^2\tau^2\|\bX^{(t+1)}-\bX^{(t)}\|^2_F
\\\nonumber
&+3\|\bX^{(t)}(\bY^{(t)})^{\T}-\bZ\|^2_F\|\bY^{(t+1)}-\bY^{(t)}\|^2_F
\\
&+3N\tau\|\bX^{(t)}(\bY^{(t+1)}-\bY^{(t)})^{\T}\|^2_F.
\end{align}
\end{lemma}
\begin{IEEEproof}
The optimality condition of the $\bX$ subproblem \eqref{eq.pupdatex} is given
by
\begin{equation}
\def\lhs{ (\bX^{(t+1)}(\bY^{(t+1)})^{\T}-\bZ)\bY^{(t+1)} }
\def\rhs{ +\rho(\bX^{(t+1)}-\bY^{(t+1)}+\bLambda^{(t)}/\rho)=0. }
\iftoggle{twocolumn}{
  \begin{split} & \lhs\\ & \qquad \rhs \end{split}
}{
  \lhs \rhs
}
\label{eq:opt3b}
\end{equation}
Substituting \eqref{eq.pupdatelam} into \eqref{eq:opt3b}, we have
\begin{equation}\label{eq:optlam}
\bLambda^{(t+1)}=-(\bX^{(t+1)}(\bY^{(t+1)})^{\T}-\bZ)\bY^{(t+1)}.
\end{equation}
Subtracting the same equation in iteration $t$, we have the successive
difference of the dual matrix \eqref{eq:difflambda}, shown at the top of the
next page.
\begin{figure*}[htb]
\rule{\textwidth}{.5pt}
\begin{align}
&\bLambda^{(t+1)}-\bLambda^{(t)}=-\left[\bX^{(t+1)}(\bY^{(t+1)})^{\T}\bY^{(t+1)}-\bX^{(t)}(\bY^{(t)})^{\T}\bY^{(t)}-\bZ(\bY^{(t+1)}-\bY^{(t)})\right]\label{eq:difflam}
\\\nonumber
=&-\left[(\bX^{(t+1)}-\bX^{(t)})(\bY^{(t+1)})^{\T}\bY^{(t+1)}+\bX^{(t)}\big((\bY^{(t+1)})^{\T}\bY^{(t+1)}-(\bY^{(t)})^{\T}\bY^{(t)}\big)+\bZ(\bY^{(t+1)}-\bY^{(t)})\right]
\\\nonumber
=&\bZ(\bY^{(t+1)}-\bY^{(t)})-(\bX^{(t+1)}-\bX^{(t)})(\bY^{(t+1)})^{\T}\bY^{(t+1)}\label{eq:difflambda}
\\
&-\underbrace{\frac{1}{2}\big(\bX^{(t)}\big((\bY^{(t+1)}+\bY^{(t)})^{\T}(\bY^{(t+1)}-\bY^{(t)})+(\bY^{(t+1)}-\bY^{(t)})^{\T}(\bY^{(t+1)}+\bY^{(t)})\big)\big)}_{\bydef\mathcal{Q}}.
\end{align}
\rule{\textwidth}{.5pt}
\end{figure*}

Note that the following is true
\begin{align}
\nonumber
\mathcal{Q}=&\frac{1}{2}\big(\bX^{(t)}(\bY^{(t+1)}-\bY^{(t)})^{\T}(\bY^{(t+1)}-\bY^{(t)})
\\\nonumber
&+2\bX^{(t)}(\bY^{(t)})^{\T}(\bY^{(t+1)}-\bY^{(t)})\big)
\\\nonumber
&+\frac{1}{2}\bX^{(t)}(\bY^{(t+1)}-\bY^{(t)})^{\T}(\bY^{(t+1)}+\bY^{(t)})
\\\nonumber
=&\bX^{(t)}(\bY^{(t)})^{\T}(\bY^{(t+1)}-\bY^{(t)})
\\
&+\bX^{(t)}(\bY^{(t+1)}-\bY^{(t)})^{\T}\bY^{(t+1)}.\label{eq:qexpression}
\end{align}
Plugging \eqref{eq:qexpression} into \eqref{eq:difflambda}, we have
\begin{align}
\onlytwo{&}\nonumber
\bLambda^{(t+1)}-\bLambda^{(t)}
\onlytwo{\\\nonumber}
=&\bZ(\bY^{(t+1)}-\bY^{(t)})-(\bX^{(t+1)}-\bX^{(t)})(\bY^{(t+1)})^{\T}\bY^{(t+1)}
\\\nonumber
&-\bX^{(t)}(\bY^{(t)})^{\T}(\bY^{(t+1)}-\bY^{(t)})
\onlytwo{\\\notag&}
-\bX^{(t)}(\bY^{(t+1)}-\bY^{(t)})\bY^{(t+1)}
\\\nonumber
=&(\bZ-\bX^{(t)}(\bY^{(t)})^{\T})(\bY^{(t+1)}-\bY^{(t)})
\onlytwo{\\\nonumber&}
-(\bX^{(t+1)}-\bX^{(t)})(\bY^{(t+1)})^{\T}\bY^{(t+1)}
\\
&-\bX^{(t)}(\bY^{(t+1)}-\bY^{(t)})^{\T}\bY^{(t+1)}.
\end{align}
Using triangle inequality, we arrive at
\begin{align}\label{eq:difflambdan}
\nonumber
\|\bLambda^{(t+1)}-&\bLambda^{(t)}\|_F\le\|\bX^{(t+1)}-\bX^{(t)}\|_F\|(\bY^{(t+1)})^{\T}\bY^{(t+1)}\|_F
\\\nonumber
&+\|\bX^{(t)}(\bY^{(t)})^{\T}-\bZ\|_F\|\bY^{(t+1)}-\bY^{(t)}\|_F
\\
&+\|\bX^{(t)}(\bY^{(t+1)}-\bY^{(t)})^{\T}\|_F\|\bY^{(t+1)}\|_F.
\end{align}
Since $\|\bY_i\|^2\le\tau$, we know that $\|\bY\|_F\le \sqrt{N\tau}$. Squaring
both sides of \eqref{eq:difflambdan}, we obtain
\begin{align}
\nonumber
\|\bLambda^{(t+1)}-&\bLambda^{(t)}\|^2_F\le3N^2\tau^2\|\bX^{(t+1)}-\bX^{(t)}\|^2_F
\\\nonumber
&+3\|\bX^{(t)}(\bY^{(t)})^{\T}-\bZ\|^2_F\|\bY^{(t+1)}-\bY^{(t)}\|^2_F
\onlytwo{\\&}
+3N\tau\|\bX^{(t)}(\bY^{(t+1)}-\bY^{(t)})^{\T}\|^2_F.
\end{align}
The claim is proved.
\end{IEEEproof}

In the second step, we bound the successive difference of the augmented
Lagrangian.
\begin{lemma}\label{le.nondecrease}
Consider using the update rules \eqref{eq.pupdatey}--\eqref{eq.pupdatelam}. If
\begin{equation}\label{eq:optcondition}
\rho>
6N\tau\quad\textrm{and}\quad\beta^{(t)}>\frac{6}{\rho}\|\bX^{(t)}(\bY^{(t)})^{\T}-\bZ\|^2_F-\rho,
\end{equation}
we have
\begin{equation}
\def\parta{\mathcal{L}(\bX^{(t+1)},\bY^{(t+1)},\bLambda^{(t+1)})
  -\mathcal{L}(\bX^{(t)},\bY^{(t)},\bLambda^{(t)})}
\def\partb{\le-c_1\|\bX^{(t+1)}-\bX^{(t)}\|^2_F-c_2\|\bX^{(t)}(\bY^{(t+1)}
  -\bY^{(t)})^{\T}\|^2_F}
\def\partc{-c_3\|\bY^{(t+1)}-\bY^{(t)}\|^2_F}
\begin{split}
  & \parta \\ & \quad \partb \onlytwo{\\ & \quad} \partc
\end{split}
\label{eq:ledesabc}
\end{equation}
where $c_1,c_2,c_3>0$ are some positive constants.
\end{lemma}
\begin{IEEEproof}
Let
\begin{equation}
\def\parta{
  \widehat{\mathcal{L}}(\bX^{(t)},\bY,\bLambda^{(t)})
    \bydef \frac{1}{2}\|\bX^{(t)}\bY^{\T}-\bZ\|^2_F}
\def\partb{+\frac{\rho}{2}\|\bX^{(t)}-\bY+\bLambda^{(t)}/\rho\|^2_F
  +\frac{\beta^{(t)}}{2}\|\bY-\bY^{(t)}\|^2_F,}
\iftoggle{twocolumn}{
  \begin{split}
  & \parta \\ &\quad \partb
  \end{split}
}{
  \parta \partb
}
\end{equation}
which is an upper bound of $\mathcal{L}(\bX^{(t)},\bY,\bLambda^{(t)})$, and
\begin{align*}
\cA &\bydef \mathcal{L}(\bX^{(t)},\bY^{(t+1)},\bLambda^{(t)})
  -\mathcal{L}(\bX^{(t)},\bY^{(t)},\bLambda^{(t)}), \\
\cB & \bydef \mathcal{L}(\bX^{(t+1)},\bY^{(t+1)},\bLambda^{(t)})
  -\mathcal{L}(\bX^{(t)},\bY^{(t+1)},\bLambda^{(t)}), \\
\cC & \bydef \mathcal{L}(\bX^{(t+1)},\bY^{(t+1)},\bLambda^{(t+1)})
  -\mathcal{L}(\bX^{(t+1)},\bY^{(t+1)},\bLambda^{(t)}), \\
\hcA & \bydef \widehat{\mathcal{L}}(\bX^{(t)},\bY^{(t+1)},\bLambda^{(t)})
  -\mathcal{L}(\bX^{(t)},\bY^{(t)},\bLambda^{(t)}).
\end{align*}
We have the following descent estimate
\begin{multline}
\mathcal{L}(\bX^{(t+1)},\bY^{(t+1)},\bLambda^{(t+1)})-\mathcal{L}(\bX^{(t)},\bY^{(t)},\bLambda^{(t)}) \\
={ \cA + \cB + \cC} \le \widehat{ \cA} + \cB + \cC.
\label{eq:desla}
\end{multline}

Next
we bound the quantities in \eqref{eq:desla}
\begin{align}
\nonumber
\widehat{\mathcal{A}}&=\frac{1}{2}\|\bX^{(t)}(\bY^{(t+1)})^{\T}-\bZ\|^2_F-\frac{1}{2}\|\bX^{(t)}(\bY^{(t)})^{\T}-\bZ\|^2_F
\\\nonumber
&+\frac{\rho}{2}\|\bX^{(t)}-\bY^{(t+1)}+\bLambda^{(t)}/\rho\|^2_F
\\\nonumber
&-\frac{\rho}{2}\|\bX^{(t)}-\bY^{(t)}+\bLambda^{(t)}/\rho\|^2_F+\frac{\beta^{(t)}}{2}\|\bY^{(t+1)}-\bY^{(t)}\|^2_F
\\\nonumber
\mathop{=}\limits^{(a)}&\langle(\bX^{(t)}(\bY^{(t+1)})^{\T}-\bZ)\bX^{(t)},\bY^{(t+1)}-\bY^{(t)}\rangle
\\\nonumber
&-\frac{1}{2}\|\bX^{(t)}(\bY^{(t+1)}-\bY^{(t)})^{\T}\|^2_F
\\\nonumber
&+\rho\langle\bX^{(t)}-\bY^{(t+1)}+\bLambda^{(t)}/\rho,\bY^{(t+1)}-\bY^{(t)}\rangle
\\\nonumber
&-\frac{\rho}{2}\|\bY^{(t+1)}-\bY^{(t)}\|^2_F+\frac{\beta^{(t)}}{2}\|\bY^{(t+1)}-\bY^{(t)}\|^2_F
\\\nonumber
\mathop{\le}\limits^{(b)}&-\frac{1}{2}\|\bX^{(t)}(\bY^{(t+1)}-\bY^{(t)})^{\T}\|^2_F-\frac{\rho}{2}\|\bY^{(t+1)}-\bY^{(t)}\|^2_F
\\\nonumber
&-\frac{\beta^{(t)}}{2}\|\bY^{(t+1)}-\bY^{(t)}\|^2_F
\end{align}
where $(a)$ is due to the fact that Taylor expansion for quadratic problems is
exact, and $(b)$ is due to the optimality condition for problem
\eqref{eq.pupdatey}.  Similarly, we have
\begin{align}\label{eq:desbc}
\nonumber
\mathcal{B}\le&-\frac{1}{2}\|(\bX^{(t+1)}-\bX^{(t)})(\bY^{(t+1)})^{\T}\|^2_F
\\
&-\frac{\rho}{2}\|\bX^{(t+1)}-\bX^{(t)}\|^2_F,
\\\nonumber
\mathcal{C}=&\langle\bX^{(t+1)}-\bY^{(t+1)},\bLambda^{(t+1)}-\bLambda^{(t)}\rangle
\\
\mathop{=}\limits^{(a)}&\frac{1}{\rho}\|\bLambda^{(t+1)}-\bLambda^{(t)}\|^2_F \label{eq.desc}
\end{align}
where $(a)$ is from \eqref{eq.pupdatelam}.

Substituting the result of \leref{le.bdofdual} into \eqref{eq.desc}, we can
obtain
\begin{align}\label{eq:desabc}
\nonumber
&\mathcal{L}(\bX^{(t+1)},\bY^{(t+1)},\bLambda^{(t+1)})-\mathcal{L}(\bX^{(t)},\bY^{(t)},\bLambda^{(t)})
\\\nonumber
&\le-\left(\frac{\rho}{2}-\frac{3N^2\tau^2}{\rho}\right)\|\bX^{(t+1)}-\bX^{(t)}\|^2_F
\\\nonumber
&-\left(\frac{1}{2}-\frac{3N\tau}{\rho}\right)\|\bX^{(t)}(\bY^{(t+1)}-\bY^{(t)})^{\T}\|^2_F
\\\nonumber
&-\left(\frac{\rho}{2}+\frac{\beta^{(t)}}{2}-\frac{3\|\bX^{(t)}(\bY^{(t)})^{\T}-\bZ\|^2_F}{\rho}\right)\|\bY^{(t+1)}-\bY^{(t)}\|^2_F
\\
&-\frac{1}{2}\|(\bX^{(t+1)}-\bX^{(t)})(\bY^{(t+1)})^{\T}\|^2_F.
\end{align}

Therefore, from \eqref{eq:desabc} if
$\frac{\rho}{2}-\frac{3N^2\tau^2}{\rho}>0$,
$\frac{1}{2}-\frac{3N\tau}{\rho}>0$, and
\begin{equation}
\frac{\rho+\beta^{(t)}}{2}-\frac{3\|\bX^{(t)}(\bY^{(t)})^{\T}-\bZ\|^2_F}{\rho}>0,
\end{equation}
which are equivalent to
\begin{equation}
\rho>
6N\tau\quad\textrm{and}\quad\beta^{(t)}>\frac{6\|\bX^{(t)}(\bY^{(t)})^{\T}-\bZ\|^2_F-\rho^2}{\rho},
\end{equation}
then
$\mathcal{L}(\bX^{(t+1)},\bY^{(t+1)},\bLambda^{(t+1)})-\mathcal{L}(\bX^{(t)},\bY^{(t)},\bLambda^{(t)})<0$.

Then, it is concluded that
$\mathcal{L}(\bX^{(t+1)},\bY^{(t+1)},\bLambda^{(t+1)})$ is decreasing.
\end{IEEEproof}

In the next step we prove that
$\mathcal{L}(\bX^{(t+1)},\bY^{(t+1)},\bLambda^{(t+1)})$ is lower bounded.

\begin{lemma}\label{le.lowerbd}
Consider using the update rules \eqref{eq.pupdatey} \eqref{eq.pupdatex}
\eqref{eq.pupdatelam}. If $\rho\ge N\tau$ is satisfied, we have
\begin{equation}\label{eq:lowerbd}
\mathcal{L}(\bX^{(t+1)},\bY^{(t+1)},\bLambda^{(t+1)})\ge 0.
\end{equation}
\end{lemma}

\begin{IEEEproof} At iteration $t+1$, the augmented Lagrangian can be lower bounded as
\begin{align}
\nonumber
&\mathcal{L}(\bX^{(t+1)},\bY^{(t+1)},\bLambda^{(t+1)})
\\\nonumber
=&\frac{1}{2}\|\bX^{(t+1)}(\bY^{(t+1)})^{\T}-\bZ\|^2_F+\langle\bX^{(t+1)}-\bY^{(t+1)},\bLambda^{(t+1)}\rangle
\\\nonumber
&+\frac{\rho}{2}\|\bX^{(t+1)}-\bY^{(t+1)}\|^2_F
\\\nonumber
\mathop{=}\limits^{(a)}&\frac{1}{2}\|\bX^{(t+1)}(\bY^{(t+1)})^{\T}-\bZ\|^2_F
\\\nonumber
&+\langle\bX^{(t+1)}-\bY^{(t+1)},-(\bX^{(t+1)}(\bY^{(t+1)})^{\T}-\bZ)\bY^{(t+1)}\rangle
\\\nonumber
&+\frac{\rho}{2}\|\bX^{(t+1)}-\bY^{(t+1)}\|^2_F
\\
\mathop{\ge}\limits^{(b)}&\frac{1}{2}(\rho-N\tau)\|\bX^{(t+1)}-\bY^{(t+1)}\|^2_F \label{eq:lowerbdcon}
\end{align}
where $(a)$ is due to \eqref{eq:optlam}, and $(b)$ is true because
\begin{align}
\nonumber
0\le&\|(\bX^{(t+1)}-\bY^{(t+1)})(\bY^{(t+1)})^{\T}-(\bX^{(t+1)}(\bY^{(t+1)})^{\T}-\bZ)\|^2_F
\\\nonumber
=&\|(\bX^{(t+1)}-\bY^{(t+1)})(\bY^{(t+1)})^{\T}\|^2_F
\\\nonumber
&-2\langle(\bY^{(t+1)})^{\T}(\bX^{(t+1)}-\bY^{(t+1)}),\bX^{(t+1)}(\bY^{(t+1)})^{\T}-\bZ\rangle
\\\nonumber
&+\|\bX^{(t+1)}(\bY^{(t+1)})^{\T}-\bZ)\|^2_F,
\end{align}
{and $\|\bY\|^2_F\le N\tau$}.

From \eqref{eq:lowerbdcon}, we know that if $\rho\ge N\tau$, we have
$\mathcal{L}(\bX^{(t+1)},\bY^{(t+1)},\bLambda^{(t+1)})\ge 0$.
\end{IEEEproof}

These lemmas lead to the main convergence claim.

\begin{IEEEproof}
Combing \eqref{eq:ledesabc} and \eqref{eq:lowerbd}, we have
\begin{align}\label{eq:limitsxy}
&\lim_{t\to\infty}\|\bX^{(t+1)}-\bX^{(t)}\|^2_F=0,
\\\nonumber
&\lim_{t\to\infty}\|\bX^{(t)}(\bY^{(t+1)}-\bY^{(t)})^{\T}\|^2_F=0,
\\\nonumber
&\lim_{t\to\infty}\|\bX^{(t)}(\bY^{(t)})^{\T}-\bZ\|^2_F\|\bY^{(t+1)}-\bY^{(t)}\|^2_F=0.
\end{align}

By \leref{le.bdofdual}, we have
\begin{equation}
\lim_{t\to\infty} \|\bLambda^{(t+1)}-\bLambda^{(t)}\|^2_F =0,
\end{equation}
which implies $\lim_{t\to\infty}\|\bX^{(t)}-\bY^{(t)}\|^2_F=0$. Combining with
\eqref{eq:limitsxy}, we can further know that
$\lim_{t\to\infty}\|\bY^{(t+1)}-\bY^{(t)}\|^2_F=0$. The boundedness assumption
of $\mathbf{X}^{(t)}$ then follows from the boundedness of $\mathbf{Y}^{(t)}$.
Using the expression of $\bLambda^{(t)}$ in \eqref{eq:optlam}, one can show
that $\{\bLambda^{(t)}\}$ is also bounded.

The optimality condition of \eqref{eq.pupdatey} is given by
\begin{multline}
\big\langle(\bX^{(t)})^{\T}(\bX^{(t)}(\bY^{(t+1)})^{\T}-\bZ)-\rho(\bX^{(t)}-\bY^{(t+1)}+\bLambda^{(t)}/\rho)^{\T}
\\ \quad
+\beta^{(t)}(\bY^{(t+1)}-\bY^{(t)})^{\T},(\bY-\bY^{(t+1)})^{\T}\big\rangle\ge
0,\\
\qquad\qquad\forall\;\bY\ge0\quad\textrm{and}\quad\|\bY_i\|^2_2\le\tau\;\forall
i.
\label{eq:optofy}
\end{multline}

Substituting \eqref{eq:optlam} into \eqref{eq:optofy}, using
\eqref{eq:limitsxy}, and taking limit over any converging subsequence of
$\{\bX^{(t)}, \bY^{(t)}, \bLambda^{(t)}\}$, we have
\begin{equation}
\begin{split}
\langle(\bX^*)^{\T}(\bX^*(\bY^*)^{\T}&-\bZ)+((\bX^*(\bY^*)^{\T}-\bZ)\bY^*)^{\T} \onlytwo{\\&}
  -\rho(\bX^*-\bY^*)^{\T},(\bY-\bY^*)^{\T}\rangle\ge0, \\
& \quad\forall\;\bY\ge0\quad\textrm{and}\quad\|\bY_i\|^2_2\le\tau\;\forall i.
\end{split}
\end{equation}

The optimality condition of \eqref{eq.pupdatex} is given by
\begin{align}
(\bX^{(t+1)}(\bY^{(t+1)})^{\T}-\bZ)(\bY^{(t+1)})
\onlytwo{\hspace{20pt}\notag \\}+\rho(\bX^{(t+1)}-\bY^{(t+1)}+\bLambda^{(t)}/\rho)=0.
\label{eq:kktx}
\end{align}
Taking limit of \eqref{eq:kktx} over the same subsequence, we have
\begin{equation}
(\bX^*(\bY^*)^{\T}-\bZ)\bY^*+\rho(\bX^*-\bY^*+\bLambda^*/\rho)=0.
\end{equation}
Using the fact $\bX^*=\bY^*$, we have
\begin{align}
\nonumber
&\bigg\langle\big(\bX^*(\bX^*)^{\T}-\frac{\bZ^{\T}+\bZ}{2}\big)\bX^*,\bX-\bX^*\bigg\rangle\ge0,
\\
&\qquad\qquad\qquad\qquad\qquad\forall\;\bX\ge0,\;\|\bX_i\|^2_2\le\tau\;\forall i,
\\
&(\bX^*(\bX^*)^{\T}-\bZ)\bX^*+\bLambda^*=0,
\end{align}
which are the KKT conditions of problem \eqref{eq:symnmf}.
\end{IEEEproof}

\subsection{Convergence Rate Proof of the Proposed Algorithm}
\label{sec:conrateproof}

\begin{IEEEproof}
Based on Theorem~\ref{th.convergence}, $\|\bX^{(t)}\|^2_F$ is bounded. There
must exist a finite $\gamma>0$ such that $\|\bX^{(t)}\|^2_F\le N\gamma,
\forall t$, where $\gamma$ is only dependent on $\tau$, $N$ and $\|\bZ\|_F$.

From the optimality condition of $\bY$ in \eqref{eq.pupdatey}, we have
\begin{equation}\notag
\begin{split}
(\bY^{(t+1)})^{\T}&=\textsf{proj}_{\mathcal{Y}}\bigg[(\bY^{(t+1)})^{\T}
  \onlytwo{\\ &}
  -((\bX^{(t)})^{\T}(\bX^{(t)}(\bY^{(t+1)})^{\T}-\bZ)-\rho(\bX^{(t)} \\
  &\;-\bY^{(t+1)}+\bLambda^{(t)}/\rho)^{\T}
  +\beta^{(t)}(\bY^{(t+1)}-\bY^{(t)})^{\T})\bigg].
\end{split}
\end{equation}
Then, we have
\begin{align}
\nonumber
&\bigg\|(\bY^{(t)})^{\T}-\textsf{proj}_{\mathcal{Y}}\big[(\bY^{(t)})^{\T}-((\bX^{(t)})^{\T}(\bX^{(t)}(\bY^{(t)})^{\T}-\bZ)
\onlytwo{\\\nonumber&}
-\rho(\bX^{(t)}-\bY^{(t)}+\bLambda^{(t)}/\rho)^{\T})\big]\bigg\|_F
\\\nonumber
&=\bigg\|(\bY^{(t)})^{\T}-(\bY^{(t+1)})^{\T}+(\bY^{(t+1)})^{\T}
\\\nonumber&\quad
-\textsf{proj}_{\mathcal{Y}}\big[(\bY^{(t)})^{\T}-((\bX^{(t)})^{\T}(\bX^{(t)}(\bY^{(t)})^{\T}-\bZ)
\onlytwo{\\\nonumber&}
\quad-\rho(\bX^{(t)}-\bY^{(t)}+\bLambda^{(t)}/\rho)^{\T})\big]\bigg\|_F
\\\nonumber
&\mathop{\le}\limits^{(a)}\|\bY^{(t)}-\bY^{(t+1)}\|_F
\\\nonumber
&\quad +\bigg\|\textsf{proj}_{\mathcal{Y}}\big[(\bY^{(t+1)})^{\T}-((\bX^{(t)})^{\T}(\bX^{(t)}(\bY^{(t+1)})^{\T}-\bZ)
\\\nonumber
&\quad -\rho(\bX^{(t)}-\bY^{(t+1)}+\bLambda^{(t)}/\rho)^{\T}+\beta^{(t)}(\bY^{(t+1)}-\bY^{(t)})^{\T})\big]
\\\nonumber
&\quad -\textsf{proj}_{\mathcal{Y}}\big[(\bY^{(t)})^{\T}-((\bX^{(t)})^{\T}(\bX^{(t)}(\bY^{(t)})^{\T}-\bZ)
\onlytwo{\\\nonumber&\quad}
-\rho(\bX^{(t)}-\bY^{(t)}+\bLambda^{(t)}/\rho)^{\T})\big]\bigg\|_F
\\\nonumber
& \mathop{\le}\limits^{(b)}(2+\rho+\beta^{(t)})\|\bY^{(t+1)}-\bY^{(t)}\|_F
\onlytwo{\\\nonumber &\quad}
+\|(\bX^{(t)})^{\T}\bX^{(t)}(\bY^{(t+1)}-\bY^{(t)})^{\T}\|_F
\\\nonumber
&\mathop{\le}\limits^{(c)}(2+\rho+\beta^{(t)})\|\bY^{(t+1)}-\bY^{(t)}\|_F
\onlytwo{\\&\quad}
+\sqrt{N\gamma}\|\bX^{(t)}(\bY^{(t+1)}-\bY^{(t)})^{\T}\|_F
\label{eq.nablaybd}
\end{align}
where $\textsf{proj}_{\mathcal{Y}}$ denotes the projection of $\bY$ to the
feasible space; in $(a)$ we used triangle inequality; $(b)$ is due to the
nonexpansiveness of the projection operator; and $(c)$ is due to the
boundedness of $\|\bX\|_F$.

Similarly, we can bound the size of the gradient of the augmented Lagrangian
with respect to $\bX$ by the following series of inequalities
\begin{align}
\nonumber
\onlytwo{&}\|\nabla_{\bX}\mathcal{L}(\bX^{(t)},\bY^{(t)},\bLambda^{(t)})\|_F
\onlyone{&}=\|(\bX^{(t)}(\bY^{(t)})^{\T}-\bZ)\bY^{(t)}
\onlytwo{\\\nonumber&}
+\rho(\bX^{(t)}-\bY^{(t)}+\bLambda^{(t)}/\rho)\|_F
\\\nonumber
&\stackrel{(a)}=\big\|(\bX^{(t)}(\bY^{(t)})^{\T}-\bZ)\bY^{(t)}+\rho(\bX^{(t)}-\bY^{(t)}+\bLambda^{(t)}/\rho)
\\\nonumber
&\quad -((\bX^{(t+1)}(\bY^{(t+1)})^{\T}-\bZ)\bY^{(t+1)}
\onlytwo{\\\nonumber&}
\quad+\rho(\bX^{(t+1)}-\bY^{(t+1)}+\bLambda^{(t)}/\rho))\big\|_F\notag
\end{align}
\begin{align}
&\le\|(\bX^{(t)}(\bY^{(t)})^{\T}-\bZ)\bY^{(t)}
\onlytwo{\notag \\ &\quad}
-((\bX^{(t+1)}(\bY^{(t+1)})^{\T}-\bZ)\bY^{(t+1)})\|_F \notag \\
&\quad +\rho\|\bY^{(t+1)}-\bY^{(t)}\|_F+\rho\|\bX^{(t+1)}-\bX^{(t)}\|_F
\\\nonumber
& \mathop{=}\limits^{(b)}\|\bLambda^{(t+1)}-\bLambda^{(t)}\|_F+\rho\|\bY^{(t+1)}-\bY^{(t)}\|_F
\\
&\quad\quad+\rho\|\bX^{(t+1)}-\bX^{(t)}\|_F \label{eq.ratebd}
\end{align}
where $(a)$ is from the optimality condition of the $\bX$-subproblem
\eqref{eq:opt3b}; $(b)$ is true due to \eqref{eq:difflam} and
\eqref{eq:optlam}. Squaring both sides of \eqref{eq.ratebd} and applying
\leref{le.bdofdual}, we have
\begin{equation}\label{eq.nablaxbd}
\begin{split}
\!\!\!&\|\nabla_{\bX}\mathcal{L}(\bX^{(t)},\bY^{(t)},\bLambda^{(t)})\|^2_F
\onlytwo{\\&}
  \le 3(3N^2\tau^2+\rho^2)\|\bX^{(t+1)}-\bX^{(t)}\|^2_F \\
  & +3(3\|\bX^{(t)}(\bY^{(t)})^{\T}-\bZ\|^2_F+\rho^2)\|\bY^{(t+1)}
    -\bY^{(t)}\|^2_F\\
  & +9N\tau\|\bX^{(t)}(\bY^{(t+1)}-\bY^{(t)})^{\T}\|^2_F.
\end{split}
\end{equation}

Due to the boundedness of $\mathbf{X}^{(t)}$ and $\bY^{(t)}$, we must have
that for some $\delta>0$, $\|\bX^{(t)}(\bY^{(t)})^{\T}-\bZ\|_F\le \delta$.

Therefore, combining \eqref{eq.nablaybd} and \eqref{eq.nablaxbd}, there must
exists a finite positive number $\sigma_1$ such that
\begin{equation}\label{eq.objpbd}
\|\widetilde{\nabla}\mathcal{L}(\bX^{(t)},\bY^{(t)},\bLambda^{(t)})\|^2_F\le\sigma_1\mathcal{F}
\end{equation}
where
\begin{equation}
\begin{split}
\mathcal{F}\bydef\|\bX^{(t+1)}-\bX^{(t)}\|^2_F+\|\bY^{(t+1)}
  -\bY^{(t)}\|^2_F \onlytwo{\\}
+\|\bX^{(t)}(\bY^{(t+1)}-\bY^{(t)})^{\T}\|^2_F
\end{split}
\end{equation}
In particular, we have $\sigma_1\bydef\max\{3(3N^2\tau^2+\rho^2),
3(2+\rho+\beta^{(t)})^2+3(3\delta^2+\rho^2),3\gamma+9N\tau\}$ and
$\beta^{(t)}\le 6\delta^2/\rho$.

According to \leref{le.bdofdual}, we have
\begin{equation}
\|\bX^{(t+1)}-\bY^{(t+1)}\|^2_F=\frac{1}{\rho^2}\|\bLambda^{(t+1)}-\bLambda^{(t)}\|^2_F\le\sigma_2\mathcal{F}
\end{equation}
where some constant $\sigma_2\bydef\max\{3N^2\tau^2/\rho^2,
3\delta^2/\rho^2,3N\tau/\rho^2\}$.

Also, we have
\begin{align}
\onlytwo{\nonumber &}
\|\bX^{(t)}-\bY^{(t)}\|_F \onlytwo{\\ \nonumber}
=&\|\bX^{(t)}-\bX^{(t+1)}+\bX^{(t+1)}-\bY^{(t+1)}+\bY^{(t+1)}-\bY^{(t)}\|_F
\\\nonumber
\le&\|\bX^{(t)}-\bX^{(t+1)}\|_F+\|\bX^{(t+1)}-\bY^{(t+1)}\|_F
\onlytwo{\\&}
+\|\bY^{(t+1)}-\bY^{(t)}\|_F,
\end{align}
which yields
\begin{equation}\label{eq.objdbd}
\|\bX^{(t)}-\bY^{(t)}\|^2_F\le\sigma_3\mathcal{F}
\end{equation}
for $\sigma_3\bydef\max\{9N^2\tau^2/\rho^2+3,
9\delta^2/\rho^2+3,9N\tau/\rho^2\}$.

The inequalities \eqref{eq.objpbd} and \eqref{eq.objdbd} imply that
\begin{equation}\label{eq.objbd}
\|\widetilde{\nabla}\mathcal{L}(\bX^{(t)},\bY^{(t)},\bLambda^{(t)})\|^2_F+\|\bX^{(t)}-\bY^{(t)}\|^2_F\le(\sigma_1+\sigma_3)\mathcal{F}.
\end{equation}

According to \leref{le.nondecrease}, there exists a constant
$\sigma_4\bydef\min\{c_1,c_2,c_3\}$ such that
\begin{equation}\label{eq.objdebd}
\mathcal{L}(\bX^{(t)},\bY^{(t)},\bLambda^{(t)})-\mathcal{L}(\bX^{(t+1)},\bY^{(t+1)},\bLambda^{(t+1)})\ge\sigma_4\mathcal{F}.
\end{equation}

Combining \eqref{eq.objbd} and \eqref{eq.objdebd}, we have
\begin{multline}\label{eq.indiviobj}
\|\widetilde{\nabla}\mathcal{L}(\bX^{(t)},\bY^{(t)},\bLambda^{(t)}\|^2_F+\|\bX^{(t)}-\bY^{(t)}\|^2_F
\le\\\frac{\sigma_1+\sigma_3}{\sigma_4}(\mathcal{L}(\bX^{(t)},\bY^{(t)},\bLambda^{(t)})-\mathcal{L}(\bX^{(t+1)},\bY^{(t+1)},\bLambda^{(t+1)})).
\end{multline}

Summing both sides of \eqref{eq.indiviobj} over $t=1,\ldots,r$, we have
\begin{align}
\nonumber
&\sum^r_{t=1}\|\widetilde{\nabla}\mathcal{L}(\bX^{(t)},\bY^{(t)},\bLambda^{(t)})\|^2_F+\|\bX^{(t)}-\bY^{(t)}\|^2_F
\\\nonumber
\le&\frac{\sigma_1+\sigma_3}{\sigma_4}(\mathcal{L}(\bX^{(1)},\bY^{(1)},\bLambda^{(1)})-\mathcal{L}(\bX^{(t+1)},\bY^{(t+1)},\bLambda^{(t+1)}))
\\
\mathop{\le}\limits^{(a)}&\frac{\sigma_1+\sigma_3}{\sigma_4}\mathcal{L}(\bX^{(1)},\bY^{(1)},\bLambda^{(1)})
\end{align}
where $(a)$ is due to \leref{le.lowerbd}.

According to the definition of $T(\epsilon)$ and
$\mathcal{P}(\bX^{(t)},\bY^{(t)},\bLambda^{(t)})$, the above inequality
becomes
\begin{equation}
T(\epsilon)\epsilon\le\frac{\sigma_1+\sigma_3}{\sigma_4}\mathcal{L}(\bX^{(1)},\bY^{(1)},\bLambda^{(1)}).
\end{equation}

Dividing both sides by $T(\epsilon)$, and by setting
$C\bydef(\sigma_1+\sigma_3)/\sigma_4$, the desired result is obtained.
\end{IEEEproof}

\subsection{Sufficient Condition of Global Optimality}\label{sec.globaloptproof}

\begin{IEEEproof}
Let $\bOmega$ be the Lagrange multipliers matrix. The Lagrangian of
problem \eqref{eq:symnmf} is given by
\begin{equation}\label{eq.lag}
\mathcal{L}(\bX,\bOmega)=\frac{1}{2}\textrm{Tr}\left((\bX\bX^{\T}-\bZ)^{\T}(\bX\bX^{\T}-\bZ)\right)-\langle\bX,\bOmega\rangle.
\end{equation}

Let $(\bX^*,\bOmega^*)$ be a KKT point of problem \eqref{eq:symnmf}. To show
global optimality of $(\bX^*, \bOmega^*)$, it is sufficient to prove the
following saddle point condition \cite[pp.~238]{bova04}
\begin{equation}\label{eq.goalofopt}
\mathcal{L}(\bX^*,\bOmega)\le\mathcal{L}(\bX^*,\bOmega^*)\le\mathcal{L}(\bX,\bOmega^*), \; \forall~\bOmega\ge 0, \;\forall~\bX.
\end{equation}

To show the left hand side of \eqref{eq.goalofopt}, we have the following
\begin{equation}
\begin{split}
\mathcal{L}(\bX^*,\bOmega^*)-\mathcal{L}(\bX^*,\bOmega)
=-\langle\bX^*,\bOmega^*\rangle-(-\langle\bX^*,\bOmega\rangle)
\onlytwo{\\}
=\langle\bX^*,\bOmega-\bOmega^*\rangle\mathop{=}\limits^{(a)}\langle\bX^*,\bOmega\rangle\mathop{\ge}\limits^{(b)} 0.
\end{split}
\end{equation}
where $(a)$ is due to \eqref{eq.kktslack}, and $(b)$ is due to $\bOmega\ge 0$
and \eqref{eq.kktoptx}.

Next we show the right hand side of \eqref{eq.goalofopt}
\begin{align}
\nonumber
\onlytwo{&}\mathcal{L}(\bX,\bOmega^*)-\mathcal{L}(\bX^*,\bOmega^*)
\onlytwo{\\\nonumber}
=&\underbrace{\frac{1}{2}\textrm{Tr}[(\bX\bX^{\T}-\bX^*(\bX^*)^{\T})(\bX\bX^{\T}-\bX^*(\bX^*)^{\T})]}_{\bydef\mathcal{M}}
\\\nonumber
&+\textrm{Tr}[(\bX^*(\bX^*)^{\T}-\bZ^{\T})(\bX\bX^{\T}-\bX^*(\bX^*)^{\T})]
\\\nonumber
&-\langle\bX-\bX^*,\bOmega^*\rangle\label{eq.kktinequa}\\
\mathop{\ge}\limits^{(a)}&\langle\bX-\bX^*,\left(\bX^*(\bX^*)^{\T}-\frac{\bZ^{\T}+\bZ}{2}\right)(\bX+\bX^*)\rangle
\onlytwo{\\\nonumber&}
-\langle\bX-\bX^*,\bOmega^*\rangle
\\\nonumber
\mathop{=}\limits^{(b)}&\langle\bX-\bX^*,\left(\bX^*(\bX^*)^{\T}-\frac{\bZ^{\T}+\bZ}{2}\right)(\bX-\bX^*)\rangle
\\
=&\textrm{Tr}\bigl[(\bX-\bX^*)^{\T}\underbrace{\left(\bX^*(\bX^*)^{\T}-\frac{\bZ^{\T}+\bZ}{2}\right)}_{\bydef\mathbf{S}}(\bX-\bX^*)\bigr]\label{eq.globaloptcond}
\end{align}
where $(a)$ is due to $\mathcal{M}\ge0$ and the fact that
\begin{equation}
\begin{split}
\bX\bX^{\T}-\bX^*(\bX^*)^{\T}=\frac{1}{2}\big[(\bX+\bX^*)(\bX-\bX^*)^{\T}
\onlytwo{\\}
+(\bX-\bX^*)(\bX+\bX^*)^{\T}\big];
\end{split}
\end{equation}
$(b)$ is true because of \eqref{eq.kktgrad}. Clearly, if we have
$\mathbf{S}\succeq0$, then the following inequality must be true
$$\mathcal{L}(\bX,\bOmega^*)-\mathcal{L}(\bX^*,\bOmega^*)\ge0.$$ This
completes the proof.
\end{IEEEproof}

\subsection{Sufficient Condition of Local Optimality}\label{sec.localoptproof}

\begin{IEEEproof}
We first simplify the term $\mathcal{M}$ in \eqref{eq.kktinequa} as follows.
\begin{align}
\nonumber
&\frac{1}{2}\textrm{Tr}[(\bX\bX^{\T}-\bX^*(\bX^*)^{\T})^{\T}(\bX\bX^{\T}-\bX^*(\bX^*)^{\T})]
\\\nonumber
\mathop{=}\limits^{(a)}&\frac{1}{2}\textrm{Tr}\bigg[\left((\bX-\bX^*)\bX^{\T}+\bX^*(\bX-\bX^*)^{\T}\right)^{\T}
\onlytwo{\\\nonumber&\qquad}
\left((\bX-\bX^*)\bX^{\T}+\bX^*(\bX-\bX^*)^{\T}\right)\bigg]
\\\nonumber
\mathop{=}\limits^{(b)}&\frac{1}{2}\textrm{Tr}\bigg[\left(\bYt(\bYt+\bX^*)^{\T}+\bX^*\bYt^{\T}\right)^{\T}
\onlytwo{\\\nonumber&\qquad}
\left(\bYt(\bYt+\bX^*)^{\T}+\bX^*\bYt^{\T}\right)\bigg]
\\\nonumber
\mathop{=}\limits^{(c)}&\frac{1}{2}\textrm{Tr}\big[\bU^{\T}\bU+\bX^*\bYt^{\T}\bU+\bYt(\bX^*)^{\T}\bU+\bX^*\bYt^{\T}\bU
\\\nonumber
&\quad+\bX^*\bYt^{\T}\bYt(\bX^*)^{\T}+\bX^*\bYt^{\T}\bX^*\bYt^{\T}
\\\nonumber
&\quad+\bYt(\bX^*)^{\T}\bU+\bYt(\bX^*)^{\T}\bYt(\bX^*)^{\T}+\bYt(\bX^*)^{\T}\bX^*\bYt^{\T}\big]
\\\nonumber
=&\frac{1}{2}\textrm{Tr}\left[\bU\bU^{\T}+4\bU\bX^*\bYt^{\T}+2\bYt(\bX^*)^{\T}\bX^*\bYt^{\T}\right]
\onlytwo{\\\nonumber&}
+\textrm{Tr}\left[\bX^*\bYt^{\T}\bX^*\bYt^{\T}\right]
\\\nonumber
=&\frac{1}{2}\textrm{Tr}\left[\bYt\left[\begin{array}{cc}\bYt^{\T} & \bI\end{array}\right]\left[\begin{array}{cc}\bI & 4\bX^*\\ \boldsymbol{0} & 2(\bX^*)^{\T}\bX^*\end{array}\right]\left[\begin{array}{cc}\bYt^{\T} & \bI\end{array}\right]^{\T}\bYt^{\T}\right]
\onlytwo{\\&\quad}
+\textrm{Tr}\left[\bX^*\bYt^{\T}\bX^*\bYt^{\T}\right]\label{eq.tightbd}
\end{align}
where $(a)$ is due to the fact that
\begin{align}
\bX\bX^{\T}-\bX^*(\bX^*)^{\T}=(\bX-\bX^*)\bX^{\T}+\bX^*(\bX-\bX^*)^{\T};
\end{align}
in $(b)$ we defined $\bYt\bydef\bX-\bX^*$ which shows the difference
between $\bX$ and $\bX^*$; and in $(c)$ we defined
$\bU\bydef\bYt\bYt^{\T}=\bU^{\T}$.

Combining \eqref{eq.globaloptcond} and \eqref{eq.tightbd}, we have
\begin{align}
\notag
\onlytwo{&}\mathcal{L}(\bX,\bOmega^*)-\mathcal{L}(\bX^*,\bOmega^*)
\onlytwo{\\\nonumber}
=&\textrm{Tr}\left[\bYt\left[\frac{1}{2}\bYt^{\T}\bYt+2\bYt^{\T}\bX^*+(\bX^*)^{\T}\bX^*\right]\bYt^{\T}\right]
\\\nonumber
&+\textrm{Tr}\left[\bX^*\bYt^{\T}\bX^*\bYt^{\T}\right]+\textrm{Tr}\left[\bYt^{\T}\left(\bX^*(\bX^*)^{\T}-\frac{\bZ^{\T}+\bZ}{2}\right)\bYt\right]
\\\nonumber
=&\sum^K_m\sum^K_n(\bYt'_m)^{\T}\mathcal{K}_{m,n}\bYt'_n+\sum^K_m\sum^K_n(\bYt'_m)^{\T}\widetilde{\mathcal{K}}_{m,n}\bYt'_n
\\\nonumber
&+\sum^K_m(\bYt'_m)^{\T}\mathbf{S}\bYt'_m
\\\nonumber
=&\textrm{vec}(\bYt)^{\T}\mathbf{T}\textrm{vec}(\bYt)
\end{align}
where
\begin{equation}\notag
\mathbf{T}\bydef\left[\begin{array}{ccc} \mathcal{K}_{1,1}\bI+\widetilde{\mathcal{K}}_{1,1}+\mathbf{S} & \cdots & \mathcal{K}_{1,K}\bI+\widetilde{\mathcal{K}}_{1,K}\\ \vdots & \cdots & \vdots\\ \mathcal{K}_{K,1}\bI+\widetilde{\mathcal{K}}_{K,1} & \cdots & \mathcal{K}_{K,K}\bI+\widetilde{\mathcal{K}}_{K,K}+\mathbf{S}
\end{array}\right],
\end{equation}
\begin{align}
\mathcal{K}_{m,n}\bydef\frac{1}{2}(\bYt'_m)^{\T}\bYt'_n+2(\bYt'_m)^{\T}\bX'^*_n+(\bX'^*_m)^{\T}\bX'^*_n,
\end{align}
and $\widetilde{\mathcal{K}}_{m,n}\bydef\bX'^*_n(\bX'^*_m)^{\T}$, $(m,n)$
denotes the $(m,n)$th block of a matrix, $\bX'^*_m$ ($\bYt'_n$) denotes the
$m$th (or $n$th) \emph{column} of matrix $\bX^*$ (or $\bYt$).

For the $(m,n)$th block, we have
\begin{align}
\nonumber
&(\bYt'_m)^{\T}\Big(\left(\frac{1}{2}(\bYt'_m)^{\T}\bYt'_n+2(\bYt'_m)^{\T}\bX'^*_n+(\bX'^*_m)^{\T}\bX'^*_n\right)\bI
\onlytwo{\\\nonumber&}
+\bX'^*_n(\bX'^*_m)^{\T}+\delta_{m,n}\mathbf{S}\Big)\bYt'_n
\\\nonumber
\mathop{\ge}\limits^{(a)}&(\bYt'_m)^{\T}\Big(\Big(-\frac{1}{4}\left(\|\bYt'_m\|^2_2+\|\bYt'_n\|^2_2\right)-\frac{1}{\delta}\|\bYt'_m\|^2_2
\onlytwo{\\\nonumber&}
-\delta\|\bX'^*_n\|^2_2+(\bX'^*_m)^{\T}\bX'^*_n\Big)\bI+\bX'^*_n(\bX'^*_m)^{\T}+\delta_{m,n}\mathbf{S}\Big)\bYt'_n
\\\nonumber
=&(\bYt'_m)^{\T}\left(-(\frac{1}{4}+\frac{1}{\delta})\|\bYt'_m\|^2_2-\frac{1}{4}\|\bYt'_n\|^2_2\right)\bYt'_n
\\\nonumber&\quad
+(\bYt'_m)^{\T}\bigg(\left((\bX'^*_m)^{\T}\bX'^*_n-\delta\|\bX'^*_n\|^2_2\right)\bI+\bX'^*_n(\bX'^*_m)^{\T}
\onlytwo{\\\nonumber&}
+\delta_{m,n}\mathbf{S}\bigg)\bYt'_n
\\\nonumber
\mathop{\ge}\limits^{(b)}&\|\bYt'_m\|\|\bYt'_n\|\left(-(\frac{1}{4}+\frac{1}{\delta})\|\bYt'_m\|^2_2-\frac{1}{4}\|\bYt'_n\|^2_2\right)
\onlytwo{\\\nonumber&\quad}
+(\bYt'_m)^{\T}\mathbf{T}_{m,n}\bYt'_n
\end{align}
where
\begin{equation}\notag
\mathbf{T}_{m,n}\bydef\left((\bX'^*_m)^{\T}\bX'^*_n-\delta\|\bX'^*_n\|^2_2\right)\bI+\bX'^*_n(\bX'^*_m)^{\T}+\delta_{m,n}\mathbf{S},
\end{equation}
$\delta_{m,n}$ is the Kronecker delta function, and $\mathbf{T}_{m,n}$ is the
$(m,n)$th block of matrix $\mathbf{T}$, and $(a)$ we use triangle
inequality and $\delta>0$ is any positive number; $(b)$ we use Cauchy-Schwarz
inequality.

If there exists $\delta$ such that $\mathbf{T}$ is positive definite, then
$\bX^*$ is a strict local minimum point of problem \eqref{eq:symnmf}. That is, there
exist some $\gamma, \epsilon>0$ such that
\begin{equation}\label{eq.localcm}
\begin{split}
\mathcal{L}(\bX,\bOmega^*)-\mathcal{L}(\bX^*,\bOmega^*)
  \ge\frac{\gamma}{2}\|\bX-\bX^*\|^2_F,\onlytwo{\\}
\quad\forall ~\bX~\mbox{such that}~\|\bX'_m-\bX'^*_m\|^2_2\le\epsilon,
\end{split}
\end{equation}
 where $\gamma$ is given by
\begin{equation}
\gamma=-\left(\frac{2K^2}{\delta}+K(K-2)\right)\epsilon^2+2\lambda_{\min}(\mathbf{T})
\end{equation}
where $\lambda_{\min}(\mathbf{T})$ is the smallest eigenvalue of matrix
$\mathbf{T}$. Clearly $\gamma$ can be made positive for sufficiently small
$\epsilon$.

According to the definition of Lagrangian \eqref{eq.lag}, we have
\begin{equation}
\mathcal{L}(\bX,\bOmega^*)=f(\bX)-\langle\bX,\bOmega^*\rangle.
\end{equation}

Combing with \eqref{eq.localcm} and KKT conditions
\eqref{eq.kktopto}--\eqref{eq.kktslack}, we can obtain
\begin{equation}
\begin{split}
f(\bX)\ge\mathcal{L}(\bX,\bOmega^*)
  \ge f(\bX^*)+\frac{\gamma}{2}\|\bX-\bX^*\|^2_2, \onlytwo{\\}
\quad\forall ~\bX\ge0~\mbox{such that}~ \|\bX-\bX^*\|\le\epsilon.
\end{split}
\end{equation}
Therefore $\bX^*$ is a strict local minimum point of problem
\eqref{eq:symnmf}.
\end{IEEEproof}

\subsection{Sufficient Local Optimality Condition When $K=1$}
\label{sec.localoptproofkequal1}

\begin{IEEEproof}
The term $\mathcal{M}$ is as follows.
\begin{align}\notag
\mathcal{M}=\frac{1}{2}\textrm{Tr}[\bYt[\begin{array}{cc}\bYt^{\T} &
\bI\end{array}]\left[\begin{array}{cc}\bI & 4\bX^*\\ \boldsymbol{0} &
2(\bX^*)^{\T}\bX^*\end{array}\right][\begin{array}{cc}\bYt^{\T} &
\bI\end{array}]^{\T}\bYt^{\T}]
\onlytwo{\\\hspace{2pt}}+\textrm{Tr}\left[\bX^*\bYt^{\T}\bX^*\bYt^{\T}\right]\label{eq.tightbd2}.
\end{align}

When $K=1$, \eqref{eq.tightbd2} becomes
\begin{align}
&\frac{1}{2}\byt^{\T}\byt\left[\begin{array}{cc}\byt^{\T} & 1\end{array}\right]\left[\begin{array}{cc}\bI & 4\bx^*\\ \boldsymbol{0} & 2(\bx^*)^{\T}\bx^*\end{array}\right]\left[\begin{array}{cc}\byt^{\T} & 1\end{array}\right]^{\T}
\notag \\\nonumber
&\quad+\textrm{Tr}\left[\bx^*\byt^{\T}\bx^*\byt^{\T}\right]
\\
=&\frac{1}{2}\byt^{\T}\byt\left(\byt^{\T}\byt+4\byt^{\T}\bx^*+2(\bx^*)^{\T}\bx^*\right)+\byt^{\T}\bx^*(\bx^*)^{\T}\byt
\end{align}
where $\bx^*$ and $\byt$ denote the column of matrix $\bX^*$ and $\bYt$.

Combining with \eqref{eq.globaloptcond}, we have
\begin{align}\notag
\onlytwo{&}\mathcal{L}(\bx,\bOmega^*)-\mathcal{L}(\bx^*,\bOmega^*)
\onlytwo{\\\nonumber}
&=\byt^{\T}\left[\frac{1}{2}\byt^{\T}\byt+2\byt^{\T}\bx^*+(\bx^*)^{\T}\bx^*\right]\byt
\\\nonumber
&\quad+\byt^{\T}\left[2\bx^*(\bx^*)^{\T}-\frac{\bZ^{\T}+\bZ}{2}\right]\byt
\\\nonumber
& \mathop{\ge}\limits^{(a)}\byt^{\T}\left[\frac{1}{2}\byt^{\T}\byt-\frac{1}{\delta}\|\byt\|^2_2-\delta\|\bx^*\|^2_2+(\bx^*)^{\T}\bx^*\right]\byt
\\\nonumber
&\quad+\byt^{\T}\left[2\bx^*(\bx^*)^{\T}-\frac{\bZ^{\T}+\bZ}{2}\right]\byt
\\\nonumber
&=\frac{1}{2}\|\byt\|^4_2-\frac{1}{\delta}\|\byt\|^4_2
\\\nonumber
&\quad+\byt^{\T}\underbrace{\left[\left(1-\delta\right)\|\bx^*\|^2_2\bI+2\bx^*(\bx^*)^{\T}-\frac{\bZ^{\T}+\bZ}{2}\right]}_{\bydef\mathbf{T}_1}\byt
\end{align}
where in $(a)$ we have used the triangle inequality and $\delta>0$ is any
positive number.

If there exists $\delta>0$ which ensures that $\mathbf{T}_1\succ 0$, then
there exist some $\gamma, \epsilon>0$ such that the following is true
\begin{align}
\mathcal{L}(\bx,\bOmega^*)-\mathcal{L}(\bx^*,\bOmega^*)\ge\frac{\gamma}{2}\|\bx-\bx^*\|^2_2,\onlytwo{\notag \\}
\quad\forall ~\bx~\mbox{such that}~ \|\bx-\bx^*\|\le\epsilon.\label{eq.localc}
\end{align}
In the above inequality, the constant $\gamma$ is given by
\begin{equation}
\gamma=\left(1-\frac{2}{\delta}\right)\epsilon^2+2\lambda_{\min}(\mathbf{T}_1)
\end{equation}
where $\lambda_{\min}(\mathbf{T}_1)$ denotes the smallest eigenvalue of
$\mathbf{T}_1$. Clearly $\gamma$ can be made positive by setting $\epsilon$
sufficiently small.

According to the definition of the Lagrangian, we have
\begin{equation}
\mathcal{L}(\bx,\bOmega^*)=f(\bx)-\langle\bx,\bOmega^*\rangle.
\end{equation}

Therefore, combining with \eqref{eq.localc} and the KKT conditions, we can
obtain
\begin{align}
f(\bx)\ge\mathcal{L}(\bx,\bOmega^*)\ge
f(\bx^*)+\frac{\gamma}{2}\|\bx-\bx^*\|^2_2.\onlytwo{\notag \\}
\quad\forall ~\bx\ge0~\mbox{such
that}~ \|\bx-\bx^*\|\le\epsilon.
\end{align}
\end{IEEEproof}
\small
\newpage
\bibliographystyle{IEEE-unsorted}
\bibliography{refs}
\end{document}